\documentclass[reqno]{amsart}

\numberwithin{equation}{section}

\theoremstyle{plain}
\newtheorem{theorem}{Theorem}[section]
\newtheorem{lemma}[theorem]{Lemma}
\newtheorem{corollary}[theorem]{Corollary}
\newtheorem{proposition}[theorem]{Proposition}
\newtheorem{conjecture}[theorem]{Conjecture}
\theoremstyle{definition}
\theoremstyle{remark}
\newtheorem{remark}[theorem]{Remark}
\newtheorem*{acknowledgments}{Acknowledgments}

\newcommand\nc{\newcommand}
\nc\rnc{\renewcommand}

\nc\id{\operatorname{id}}
\nc\Span{\operatorname{Span}}

\nc\modZ {{\mathbb{Z}}}
\nc\modQ {{\mathbb{Q}}}
\nc\modN {{\mathbb{N}}}
\nc\modU {{\mathcal{U}}}
\nc\modA {{\mathcal A}}
\nc\modG {{\mathcal G}}
\nc\hf{{\frac12}}

\nc\zzzcircle {\otimes}

\nc\bb[2]{\biggl[\begin{matrix}{#1}\\{#2}\end{matrix}\biggr]}
\nc\BB[2]{\{#1\}_{#2}}
\nc\BBB[2]{\{#1\}'_{#2}}

\nc\hU{{\hat \modU }}
\nc\tU{{\tilde \modU }}
\nc\dU{{\dot \modU }}
\nc\dA{{\dot \modA }}
\nc\vU{{\vec \modU }}
\nc\hA{{\hat \modA }}
\nc\cU{{\check \modU }}
\nc\bU{{\overline U}}

\nc\ho{{\widehat\otimes _h}}
\nc\plim{\varprojlim}

\nc\UA{U_\modA }
\nc\UAq{U_{\modA _q}}

\nc\floor[1]{{\lfloor#1\rfloor}}
\nc\fd[1]{\floor{\frac{#1}{d}}}
\nc\wt{\widetilde}
\nc\ta{{\wt a}}
\nc\tb{{\wt b}}
\nc\tx{{\wt x}}
\nc\ty{{\wt y}}
\nc\tz{{\wt z}}
\nc\tw{{\wt w}}

\nc\trr{\triangleright}

\nc\ver {\ |\ }
\nc\zzzcolon {\colon}
\nc\simeqto{\overset{\simeq}{\longrightarrow}}
\nc\none{(\text{none})}

\begin{document}

\title[An integral form of the quantized enveloping algebra]{An integral form of the quantized enveloping algebra of $sl_2$ and its completions}
\author{Kazuo Habiro}
\address{Research Institute for Mathematical Sciences\\ Kyoto University\\ Kyoto\\ 606-8502\\ Japan}
\email{habiro@kurims.kyoto-u.ac.jp}

\date{May 12, 2006\quad (First version: October 31, 2002)} 

\begin{abstract}
  We introduce an integral form $\modU $ of the quantized enveloping
  algebra of $sl_2$.  The algebra $\modU $ is just large enough so that
  the quasi-$R$-matrix is contained in a completion of $\modU \otimes \modU $.  We
  study several completions of the algebra $\modU $, and determine their
  centers.  This study is motivated by a study of integrality
  properties of the quantum $sl_2$ invariants of links and integral
  homology spheres.
\end{abstract}


\thanks{Partially supported by the Japan Society
  for the Promotion of Science, Grant-in-Aid for Young Scientists (B),
  16740033.}

\maketitle

\setlength{\parskip}{2pt}

\section{Introduction}
\label{sec:introduction}

The purpose of this paper is to introduce an integral form of the
quantized enveloping algebra $U_v(sl_2)$ of the Lie algebra $sl_2$,
and some completions of it.  The motivation of this paper is a study
of integrality properties of the colored Jones polynomials \cite{RT1}
of links and the Witten-Reshetikhin-Turaev invariant \cite{W,RT2} of
$3$-manifolds, which we announced in \cite{H:rims2001}.

Let $U=U_v(sl_2)$ denote the quantized enveloping algebra of the Lie
algebra $sl_2$, which is defined to be the algebra over the rational
function field $\modQ (v)$ generated by the elements $K$, $K^{-1}$, $E$,
and $F$, subject to the relations $KK^{-1}=K^{-1}K=1$ and
\[
 KE=v^2EK,\quad KF=v^{-2}FK,\quad EF-FE=\frac{K-K^{-1}}{v-v^{-1}}.
\]
Let $\modA =\modZ [v,v^{-1}]$ be the Laurent polynomial ring, and set for
$i\in \modZ $ and $n\ge 0$
\[
  [i] = \frac{v^i-v^{-i}}{v-v^{-1}}, \quad [n]! = [1][2]\cdots[n].
\]

There are at least two well-known interesting $\modZ [v,v^{-1}]$-forms of
$U$, $\UA$ and $\bU$.  $\UA$ is defined to be the $\modA $-subalgebra
$\UA$ of $U$ generated by the elements $K$, $K^{-1}$, and the divided
powers $E^{(n)}=E^n/[n]!$ and $F^{(n)}=F^n/[n]!$ for $n\ge 1$.  $\bU$ is
defined to be the $\modA $-subalgebra $\bU\subset \UA$ of $U$, which is
generated by the elements $K$, $K^{-1}$, $e=(v-v^{-1})E$, and
$f=(v-v^{-1})F$.  See \cite{Lu,D-CP} for the details of these two
$\modZ [v,v^{-1}]$-forms.  (Both $\UA$ and $\bU$ are defined for
finite-dimensional simple Lie algebras, or more generally for
Kac-Moody Lie algebras.)  These algebras inherit Hopf $\modA $-algebra
structures from a Hopf algebra structure of $U$.

Let $\Theta $ denote the quasi-$R$-matrix of $U$
\begin{equation}
  \label{eq:187}
  \Theta = \sum_{n\ge 0} (-1)^n v^{-\hf n(n-1)}(v-v^{-1})^n[n]! F^{(n)}\otimes E^{(n)},
\end{equation}
which is an element of a completion of $U\otimes U$.  For the definition of
the quasi-$R$-matrix, see \cite{Lu:book,J}.  It is well-known that one
can use $\Theta $ to obtain a left $U$-module isomorphism
\begin{equation*}
  b_{V,W}\colon V\otimes W\overset{\simeq}\longrightarrow W\otimes V
\end{equation*}
for finite dimensional left $U$-modules $V$ and $W$.  For the
category of finite dimensional left $U$-modules, the $b_{V,W}$
define a braided category structure, using which one can define the
colored Jones polynomials of links.

One of the important properties of $\Theta $ is {\em integrality} in $\UA$,
i.e., each term in the infinite sum \eqref{eq:187} is contained in
$\UA^{\otimes 2}$, and hence $\Theta $ is contained in a certain completion of
$\UA^{\otimes 2}$.  One of the important consequences of this integrality in
quantum topology is that the colored Jones polynomials of links take
values in a Laurent polynomial ring.

We introduce a yet another $\modZ [v,v^{-1}]$-form $\modU $ of $U$ which is
useful in the study of the quantum link invariants.  This algebra may
be regarded as a ``mixed version'' of $\UA$ and $\bU$.

Let $\modU $ be the $\modA $-subalgebra of $U$ generated by the elements
$K$, $K^{-1}$, $e$, and the $F^{(n)}$ for $n\ge 1$.  Obviously,
$\bU\subset \modU \subset \UA$.  We can rewrite \eqref{eq:187} as follows
\begin{equation}
  \label{eq:90}
  \Theta = \sum_{n\ge 0} (-1)^n v^{-\hf n(n-1)}F^{(n)}\otimes e^n.
\end{equation}
Observe that each term $(-1)^nv^{-\hf n(n-1)}F^{(n)}\otimes e^n$ is contained
in $\modU \otimes _\modA \modU $.  For each $n\ge 0$, let $\modU ^e_n=\modU e^n\modU $ be the two
sided ideal in $\modU $ generated by the element $e^n$.  Let $\cU$ denote
the completion
\[
  \cU = \plim_n \modU /\modU ^e_n,
\]
which will turn out to inherit from $\modU $ a complete Hopf algebra
structure over $\modA $.  The comultiplication
$\check\Delta \colon \cU\rightarrow \cU\check\otimes \cU$ takes values in the completed tensor
product
\[
  \cU\check\otimes \cU= \plim_{k,l} (\cU\otimes _\modA \cU)/
  (\overline{\modU ^e_k}\otimes _\modA \cU+\cU\otimes _\modA \overline{\modU ^e_l}),
\]
where $\overline{\modU ^e_n}$ is the closure of $\modU ^e_n$ in $\cU$.  Then
we can regard $\Theta $ as an element of $\cU\check\otimes \cU$.  Unfortunately,
the structure of $\cU$ seems quite complicated.  Therefore we also
consider some other completions of $\modU $, which have more controllable
structures than $\cU$.

For $i\in \modZ $ and $n\ge 0$, set
\[
  \{i\}= v^i-v^{-i},\quad \{n\}! = \{1\}\cdots\{n\} = (v-v^{-1})^n[n]!.
\]
Let $\hA$ and $\dA$ denote the completions of $\modA $ defined by
\begin{equation*}
  \hA = \plim_n \modA /(\{n\}!),\quad \dA = \plim_n \modA /((v-v^{-1})^n).
\end{equation*}
Since $\{n\}!$ is divisible by $(v-v^{-1})^n$, there is a natural
homomorphism $\hA\rightarrow \dA$, which is injective (see
Proposition \ref{thm:13}).  We regard $\hA\subset \dA$.  For each
$n\ge 0$, let $\modU _n$ denote the two-sided ideal in $\modU $ generated by the
elements
\begin{equation*}
  \BB{H+m}{i} e^{n-i}\quad (m\in \modZ ,\ 0\le i\le n),
\end{equation*}
where
\begin{equation*}
  \BB{H+m}{i} = \prod_{k=0}^{i-1}(v^{m-i}K-v^{-m+i}K^{-1}).
\end{equation*}
It turns out that the two-sided ideal $\modU _1$ of $\modU $ is generated by
the elements $v-v^{-1}$, $K-K^{-1}$, and $e$.  We define the
completions $\hU$ and $\dU$ of $\modU $ by
\[
  \hU = \plim_n \modU /\modU _n,\quad \dU = \plim_n \modU /(\modU _1)^n,
\]
which will be shown to have complete Hopf algebra structures over
$\hA$ and $\dA$, respectively.

Let $U_h=U_h(sl_2)$ be the $h$-adic version of $U_v(sl_2)$, for the
definition of which see Section \ref{sec:quant-group-uhsl2}.  The
algebra $U_h$ contains $\UA$, and hence $\bU$ and $\modU $.  We will show
that there is a sequence of natural homomorphisms
\[
  \cU  \rightarrow  \hU \rightarrow  \dU \rightarrow  U_h,
\]
where the two arrows $\hU \rightarrow  \dU \rightarrow  U_h$ are injective, and the arrow
$\cU\rightarrow \hU$ is conjectured to be injective.  We set
$\tU=\operatorname{Im}(\cU\rightarrow \hU)\subset \hU$.  We have
\begin{equation*}
  \tU\subset \hU\subset \dU\subset U_h.
\end{equation*}

We will determine the centers of $\modU $, $\hU$, $\dU$, and $\tU$ as
follows.  Let $C$ denote the well-known central element
\begin{equation*}
  C=(v-v^{-1})Fe+vK+v^{-1}K^{-1}\in Z(\modU ).
\end{equation*}
For $n\ge 0$, set
\[
\sigma _n=\prod_{i=1}^{n}(C^2-(v^i+v^{-i})^2).
\]
Note that $\sigma _n$ is a monic polynomial of degree $n$ in $C^2$.
\begin{theorem}
  \label{thm:27}
  The center $Z(\modU )$ of $\modU $ is freely generated as an $\modA $-algebra by
  $C$, i.e., we have $Z(\modU ) = \modA [C]$.  The centers of $\tU$, $\hU$,
  and $\dU$ are identified with completions of $\modA [C]$ as follows.
  \begin{gather*}
    Z(\tU)\simeq \plim_n \modA [C]/(\sigma _n),\quad
    Z(\hU)\simeq \plim_n \hA[C]/(\sigma _n),\\
    Z(\dU)
    \simeq \plim_n \dA[C]/(\sigma _n)
    \simeq \plim_n \dA[C]/((C^2-[2]^2)^n).
  \end{gather*}
  Therefore $Z(\tU)$ (resp. $Z(\hU)$, $Z(\dU)$) consists of the
  elements which are
  uniquely expressed as infinite sums
  \[
  z = \sum_{n\ge 0}z_n\sigma _n,
  \]
  where $z_n\in \modA +\modA C$ (resp. $z_n\in \hA+\hA C$, $z_n\in \dA+\dA C$) for
  $n\ge 0$.
\end{theorem}
The proof of Theorem \ref{thm:27} is divided into
Theorems \ref{thm:29}, \ref{thm:33}, \ref{thm:34}, \ref{thm:12}
and \ref{thm:15} below.

The results of the present paper are used in \cite{H:unified} in which
we prove integrality results of the $sl_2$ quantum invariants,
announced in \cite{H:rims2001}.  For a string knot $L$ (i.e., a
$(1,1)$-tangle consisting only of one string) with $0$ framing, let
$J_L\in Z(U_h)$ denote the {\em universal $sl_2$ invariant} of $L$
(\cite{La,Reshetikhin:89}, see also \cite{H:rims2001,H:universal}).
This invariant has a universality property over the colored Jones
polynomials of knots: The colored Jones polynomial of the knot
obtained by closing the two ends of $L$, associated to a
finite-dimensional irreducible representation $V$ of $U_h$, is
obtained from $J_L$ by taking its quantum trace in $V$.  Using
\cite[Corollary 9.15]{H:universal}, one can prove that $J_L$ is
contained in $Z(\modG _0\tU_q)$, where $\tU_q$ is a $\modZ [q,q^{-1}]$-form of
$\tU$ defined in Section \ref{sec:q-forms}, where $q=v^2$.  Using
Theorem \ref{thm:38} below, which is a slight modification of Theorem
\ref{thm:27}, we obtain the following result.  (The detailed proof
can be found in \cite{H:unified}.)

\begin{theorem}
  \label{thm:7}
  Let $L$ be a string knot with $0$ framing.  Then there are unique
  elements $a_n(L)\in  \modZ [q,q^{-1}]$ for $n\ge 0$ such that
  \[
  J_L = \sum_{n\ge 0} a_n(L) \sigma _n.
  \]
\end{theorem}

For consequences of Theorem \ref{thm:7}, see \cite{H:rims2001}.  In
particular, using Theorem \ref{thm:7} and results in the present
paper, we prove in \cite{H:unified} the existence of an invariant
$I(M)$ of integral homology $3$-spheres $M$ with values in the ring
\begin{equation*}
  \plim_n\modZ [q,q^{-1}]/((1-q)(1-q^2)\cdots(1-q^n))
\end{equation*}
which specializes at each root $\zeta $ of unity to the
Witten-Reshetikhin-Turaev invariants of $M$ at $\zeta $.

In a joint work with Le \cite{HL}, we will generalize part of this
paper to the quantized enveloping algebras of finite dimensional
semisimple Lie algebras, and prove that the invariant $I(M)$ mentioned
above generalizes to such Lie algebras.

The motivation of this work is, as explained above, from the study of
quantum invariants.  But the rest of the paper is purely algebraic.
It is organized as follows.  Section \ref{sec:quant-group-uhsl2}
recalls the definitions and the necessary properties of the quantized
enveloping algebras $U_v(sl_2)$ and $U_h(sl_2)$, as well as the
integral forms $\UA$ and $\bU$.  We will describe in Section
\ref{sec:subalgebra-w-uqsl2} the Hopf algebra structure of our
integral form $\modU $.  In Section \ref{sec:fine-hopf-algebra}, we will
introduce generalities on filtrations and completions of Hopf
algebras.  In Section \ref{sec:topology--}, we will study the
decreasing filtration $\{\modU _n\}_n$ of $\modU $, and the completion $\hU$.
In Section \ref{sec:U1-adic-topology}, we will consider the
$\modU _1$-adic filtration of $\modU $, study the completion $\dU$, and prove
the injectivity of $\hU\rightarrow \dU$ and $\dU\rightarrow U_h$.  In Section
\ref{sec:topology-c-}, we study the completion $\cU$ and the algebra
$\tU$.  Section \ref{sec:2-gradings-uq} introduces $(\modZ /2)$-gradings
on the algebra $\modU $ and the other algebras, which count ``degrees in
$K$'' in the sense that the grading for $\modU $ is as $\modU = \modG _0\modU \oplus
\modG _1\modU $, $\modG _1\modU =K\modG _0\modU $.  We determine the centers of $\modU $, $\hU$,
$\tU$, and $\dU$ in Sections \ref{sec:centers} and \ref{sec:proof-nen}
using the Harish-Chandra homomorphism.  Section \ref{sec:q-forms}
gives $\modZ [v^2,v^{-2}]$-forms of $\modU $ and its completions.

\begin{acknowledgments}
  I would like to thank Thang Le, Gregor Masbaum, Hitoshi Murakami,
  and Tomotada Ohtsuki for helpful comments and discussions.
\end{acknowledgments}

\section{The quantized enveloping algebras $U_v(sl_2)$ and $U_h(sl_2)$}
\label{sec:quant-group-uhsl2}
In this section, we will recall the definitions and properties of the
algebras $U_v(sl_2)$ and $U_h(sl_2)$, and fix notations.

\subsection{Base rings and fields}
\label{sec:base-rings-fields}
Set $q=v^2$.  We will also need the rational function field $\modQ (q)\subset \modQ (v)$
and the Laurent polynomial ring $\modA _q=\modZ [q,q^{-1}]\subset \modA $.  We have
$\modA _q=\modA \cap \modQ (q)$, and
\begin{equation*}
  \modQ (v) = \modQ (q)\oplus v\modQ (q),\quad \modA =\modA _q\oplus v\modA _q.
\end{equation*}

Let $h$ be an indeterminate, and let $\modQ [[h]]$ denote the formal power
series ring.  We will regard $\modA $ and $\modA _q$ as subrings of $\modQ [[h]]$
by setting
\[
v=\exp \frac h2 \in \modQ [[h]],\quad
q = v^2 = \exp h\in \modQ [[h]].
\]
The quotient field $\modQ ((h))$ of
$\modQ [[h]]$ contains the fields $\modQ (v)$ and $\modQ (q)$ as subfields.

\subsection{Algebra structures of $U_h$ and $U$}
\label{sec:algebra-struct-uh}

The $h$-adic quantized universal enveloping algebra $U_h=U_h(sl_2)$ is defined
to be the $h$-adically complete $\modQ [[h]]$-algebra, topologically
generated by the elements $H$, $E$, and $F$, subject to the relations
\[
 HE=E(H+2),\quad
 HF=F(H-2),\quad
 EF-FE=\frac{K-K^{-1}}{v-v^{-1}},
\]
where we set
\[
 K = v^H=\exp\hf hH.
\]
The algebras $U_h$ and $U$ are compatible in the sense that they can
naturally be regarded as subalgebras of the $\modQ ((h))$-algebra
$U_h\otimes _{\modQ [[h]]}\modQ ((h))$.

Let $U_h^0$ (resp. $U_h^+$, $U_h^-$) denote the closure of the
$\modQ [[h]]$-subalgebra of $U_h$ generated by $H$ (resp. $E$, $F$), and
let $U^0$ (resp., $U^+$, $U^-$) denote the $\modQ (v)$-subalgebra
of $U$ generated by $K^{\pm 1}$ (resp.  $E$, $F$).  We have the
triangular decompositions
\[
  U_h\simeq U_h^-\ho U_h^0 \ho U_h^+,\quad
  U\simeq U^-\otimes _{\modQ (v)} U^0 \otimes _{\modQ (v)} U^+,
\]
where $\ho$ denotes the $h$-adically completed tensor product over
$\modQ [[h]]$.

The algebra $U$ has a $\modZ $-graded $\modQ (v)$-algebra structure with the
degree determined by
\[
  \deg K^{\pm 1}=0,\quad \deg E= 1,\quad \deg F = -1.
\]
For each $n\in \modZ $, let $\Gamma _nU$ denote the degree $n$ part of $U$.
We have
\begin{equation*}
  \begin{split}
    \Gamma _nU&=\{x\in U\ver KxK^{-1}=v^{2n}x\}\\
    &=\Span_{\modQ (v)}\{F^iK^jE^k\ver j\in \modZ ,i,k\ge 0,k-i=n\}.
  \end{split}
\end{equation*}
Similar $\modZ $-grading on $U_h$ is defined and let $\Gamma _nU_h$ denote the
degree $n$ part of $U_h$.  For a homogeneous element $x$ in $U$ or
$U_h$, let $|x|$ denote the degree of $x$.  For $n\in \modZ $ and any
homogeneous additive subgroup $A$ of $U$ (resp. $U_h$), we set
\[
  \Gamma _n A = \Gamma _nU \cap  A\  (\text{resp. $\Gamma _nU_h\cap A$}).
\]

For each $j\in \modZ $, the {\em shift automorphism}
\[
  \gamma _j\colon U_h^0\rightarrow U_h^0
\]
is the unique continuous automorphism of $\modQ [[h]]$-algebra satisfying
$\gamma _j(H)=H+j$.  Similarly, the automorphism
\[
  \gamma _j\colon U^0\rightarrow U^0
\]
is defined by $\gamma _j(K)=v^jK=v^{H+j}$.  We will freely use the fact
that for any homogeneous element $x$ in $U_h$ (resp. $U$) and any
element $t$ of $U_h^0$ (resp. $U^0$), we have $t x = x\gamma _{2|x|}(t)$.

We use the following notations.  For $a\in \modZ +\modZ H$ and $n\ge 0$,
\begin{gather*}
  \{a\} = v^a-v^{-a},\\
  \BB an=\{a\}\{a-1\}\cdots \{a-n+1\},\\
  \bb an=\frac{\BB an}{\{n\}!}
  =\frac{[a][a-1]\cdots [a-n+1]}{[n]!}.
\end{gather*}
If $a=iH+j$, $i,j\in \modZ $, then  $\{iH+j\} = v^jK^i-v^{-j}K^{-i}$.

\subsection{Hopf algebra structures of $U$ and $U_h$}
\label{sec:hopf-algebra-struct-1}
The algebras $U_h$ and $U$ have compatible Hopf algebra structures
with the comultiplication $\Delta $, the counit $\epsilon $ and the antipode $S$
defined by
\begin{gather*}
    \Delta  (H) = H\otimes  1 + 1\otimes  H,\quad \epsilon (H)=0,\quad S(H)=-H,\\
    \Delta  (K) = K\otimes  K,\quad \epsilon (K)=1,\quad S(K)=K^{-1},\\
    \Delta  (E) = E\otimes  1 + K\otimes  E,\quad \epsilon (E)=0,\quad S(E)=-K^{-1}E,\\
    \Delta  (F) = F\otimes  K^{-1} + 1 \otimes  F,\quad \epsilon (F)=0,\quad S(F)=-FK.
\end{gather*}
Here $U_h$ is an $h$-adic Hopf algebra over $\modQ [[h]]$, and $U$ is an
(ordinary) Hopf algebra over $\modQ (v)$.

\subsection{The $\modA $-forms $\UA$ and $\bU$ of $U$}
\label{sec:divided-powers-e}
In the introduction, we have defined the $\modA $-subalgebras $\UA$ and
$\bU$ of $U$.  Note that they are also $\modA $-subalgebras of $U_h$.
For $*=0,+,-$, set
\begin{gather*}
  \UA^*= U^*\cap \UA=U_h^*\cap \UA,\quad  \bU^*= U^*\cap \bU=U_h^*\cap \bU.
\end{gather*}
It is known that
\begin{gather}
  \UA^0 = \Span_{\modA [K,K^{-1}]}\{\bb Hn\ver n\ge 0\},\\
  \UA^+ = \Span_\modA \{E^{(n)}\ver n\ge 0\},\quad
  \UA^- = \Span_\modA \{F^{(n)}\ver n\ge 0\},\\
  \label{eq:113}
  \UA^-\otimes _\modA  \UA^0\otimes _\modA  \UA^+ \simeqto \UA,\quad x\otimes y\otimes z\mapsto xyz,\\
  \bU^0 = \modA [K,K^{-1}],\quad  \bU^+ = \modA [e],\quad \bU^- = \modA [f],\\
  \bU^-\otimes _\modA  \bU^0\otimes _\modA  \bU^+ \simeqto \bU,\quad x\otimes y\otimes z\mapsto xyz.
\end{gather}
It is also known that $\UA$ and $\bU$ inherit Hopf $\modA $-algebra
structures from that of $U$.

\section{The $\modA $-form $\modU $ of $U$}
\label{sec:subalgebra-w-uqsl2}

\subsection{The $\modA $-form $\modU $}
\label{sec:definition-}
Recall from the introduction that $\modU $ is the $\modA $-subalgebra of $U$
generated by the elements $K$, $K^{-1}$, $e$, and the $F^{(n)}$ for
$n\ge 1$.  In other words, $\modU $ is the smallest $\modA $-subalgebra of $U$
containing $\bU^0\cup \bU^+\cup \UA^-$.  As remarked in the introduction, we
have $\bU\subset \modU \subset \UA$.  For $*=0,+,-$, set
\begin{equation*}
  \modU ^* = U^* \cap \modU  = \UA^* \cap \modU ,
\end{equation*}
which are $\modA $-subalgebras of $\modU $.

\begin{proposition}
  \label{thm:3}
  We have
  \begin{equation*}
    \modU ^0 = \bU^0,\quad
    \modU ^+ = \bU^+,\quad
    \modU ^- =\UA^-.
  \end{equation*}
  We have the triangular decomposition
  \begin{equation*}
    \modU ^-\otimes _\modA  \modU ^0\otimes _\modA  \modU ^+ \simeqto \modU ,\quad x\otimes y\otimes z\mapsto xyz.
  \end{equation*}
  Hence $\modU $ is freely $\modA $-spanned by the elements $F^{(m)}K^ie^n$
  for $i\in \modZ $, $m,n\ge 0$.
\end{proposition}

\begin{proof}
  By induction, we see that for $m,n\ge 0$
  \begin{equation}
    \label{eq:117}
    e^mF^{(n)}=\sum_{p=0}^{\min(m,n)} {\bb mp} {F^{(n-p)}}
    {\BB{H-m-n+2p}{p}} {e^{m-p}}.
  \end{equation}
  Hence $\bU^+\UA^-\subset \UA^-\bU^0\bU^+$.  We also have
  $\bU^+\bU^0=\bU^0\bU^+$ and $\bU^0\UA^-=\UA^-\bU^0$.  Then we can
  easily verify that $\UA^-\bU^0\bU^+$ is an $\modA $-subalgebra of $\modU $.
  Since $\UA^-\bU^0\bU^+$ generates $\modU $, we have
  $\UA^-\bU^0\bU^+=\modU $.  Hence there is a surjective $\modA $-module
  homomorphism $\UA^-\otimes _\modA  \bU^0\otimes _\modA  \bU^+ \rightarrow \modU $, $x\otimes y\otimes z\mapsto
  xyz$, which is injective in view of \eqref{eq:113}.  Hence
  the assertions follow.
\end{proof}

We do not need the following result in the rest
of this paper, but we state it for completeness.

\begin{proposition}
  \label{thm:2}
  As an $\modA $-algebra, $\modU $ has a presentation with generators $K$,
  $K^{-1}$, $e$, and the $F^{(n)}$ for $n\ge 1$, and with the relations
  \begin{gather}
    \label{eq:94} KK^{-1}=K^{-1}K=1,\quad
    Ke = v^2eK,\quad KF^{(n)} = v^{-2n}F^{(n)}K,\\
    \label{eq:84}
    F^{(m)}F^{(n)} = \bb{m+n}{m} F^{(m+n)}\quad (m,n\ge 1),\\
    \label{eq:82}
    eF^{(n)} = F^{(n)}e + F^{(n-1)}(v^{-n+1}K-v^{n-1}K^{-1})\quad (n\ge 1).
  \end{gather}
\end{proposition}

\begin{proof}(Sketch)
  We have the relations \eqref{eq:94}--\eqref{eq:82} in $\modU $.  We can
  prove using only these relations that $\modU $ is $\modA $-spanned by the
  elements $F^{(i)}K^je^k$ for $j\in \modZ $, $i,k\ge 0$.  It follows that the
  $\modA $-algebra $\modU $ is described by the generators and relations given
  above.
\end{proof}

\subsection{Hopf algebra structure of $\modU $}
\label{sec:hopf-algebra-struct-2}

We will use the following formulas ($m\ge 0$)
\begin{gather}
  \label{eq:h7}
  \Delta (e) = e\otimes 1+K\otimes e,\quad \epsilon (e)=0, \quad S(e)=-K^{-1}e,\\
  \label{eq:h13}
  \Delta (e^m) = \sum_{i=0}^m v^{-i(m-i)} \bb mi K^ie^{m-i}\otimes e^i,\\
  \label{eq:40}
  \epsilon (e^m)=\delta _{m,0},\quad S(e^m)= (-1)^m v^{m(m-1)} K^{-m} e^m,\\
  \label{eq:h5}
  \Delta (F^{(m)}) = \sum_{i=0}^m v^{i(m-i)} F^{(m-i)}\otimes F^{(i)}K^{-m+i},\\
  \label{eq:h8} \epsilon (F^{(m)})=\delta _{m,0},\quad
  S(F^{(m)})=(-1)^mv^{-m(m-1)}F^{(m)} K^m.
\end{gather}

Set $\modU ^{\ge 0}=\modU ^0\modU ^+ = \bU^0\bU^+$, which is an $\modA $-subalgebra of
$\modU $.

\begin{proposition}
  \label{thm:h1}
  The $\modA $-algebra $\modU $ is a Hopf $\modA $-subalgebra of $\UA$.  The
  $\modA $-subalgebras $\modU ^0$ and $\modU ^{\ge 0}$ of $\modU $ are Hopf
  $\modA $-subalgebras of $\modU $.
\end{proposition}

\begin{proof}
  It is well-known that $\modU ^0=\bU^0$ and $\modU ^{\ge 0}=\bU^0\bU^+$ are
  Hopf $\modA $-subalgebras of $\bU$.  Therefore the first assertion
  implies the second.  By \eqref{eq:h5} and \eqref{eq:h8},
  \[
  \Delta (\modU ^-)\subset \modU ^- \otimes  \modU ^-\modU ^0,\quad
  \epsilon (\modU ^-)\subset \modA ,\quad
  S(\modU ^-)\subset \modU ^-\modU ^0.
  \]
  Using the triangular decomposition, we can easily check that
  \begin{equation*}
    \Delta (\modU )\subset \modU \otimes \modU ,\quad \epsilon (\modU )=\modA ,\quad S(\modU )\subset \modU ,
  \end{equation*}
  which implies the first assertion.  The second assertion follows
  since $\modU ^0=\bU^0$ and $\modU ^{\ge 0}=\bU^0\bU^+$ are Hopf
  $\modA $-subalgebras of $\bU$.
\end{proof}

As one can easily see, $\modU ^{\le 0}=\modU ^-\modU ^0$ also is a Hopf
$\modA $-subalgebra of $\modU $, but we will not need this fact.

\section{Hopf algebra filtrations}
\label{sec:fine-hopf-algebra}

In this subsection, we fix terminology on the filtrations and
linear topology for Hopf algebras.  The contents of this section
should be well known, but we include it in the paper since we could
not find a suitable reference.

Let $R$ be a commutative ring with unit, and let
\begin{equation}
  \label{eq:F1}
  R=R_0\supset R_1\supset \cdots \supset R_n\supset \cdots
\end{equation}
be a decreasing filtration of ideals in $R$.  (Here we do {\em not}
assume that $R_mR_n\subset R_{m+n}$.)  The filtration $\{R_n\}_n$ defines a
linear topology of $R$, and the completion $\hat R=\plim_n R/R_n$ is a
commutative $R$-algebra.  If $R_1=0$, then $\{R_n\}_n$ is said to be
{\em discrete}, and we identify $\hat R$ with $R$.

Let $H=(H,\Delta ,\epsilon ,S)$ be a Hopf $R$-algebra with comultiplication $\Delta $,
counit $\epsilon $, and antipode $S$.  Suppose for simplicity that $H$ is a
free $R$-module.  A {\em Hopf algebra filtration} of $H$ with respect
to a decreasing filtration $\{R_n\}_n$ of two-sided ideals in $R$ will
mean a decreasing family of two-sided ideals in $H$
\begin{equation}
  \label{eq:F2}
  H=H_0\supset H_1\supset \cdots \supset H_n\supset \cdots
\end{equation}
such that
\begin{gather}
  \label{eq:F3} R_n\subset H_n,\\
  \label{eq:F4} \Delta (H_n) \subset  \sum_{i+j=n} H_i \otimes _R H_j,\\
  \label{eq:F5} \epsilon (H_n) \subset  R_n,\\
  \label{eq:F6} S(H_n) \subset  H_n,
\end{gather}
for each $n\ge 0$.  (Note that we do {\em not} assume
$H_mH_n\subset H_{m+n}$, $m,n\ge 0$.)  If the filtration $\{R_n\}_n$ is
discrete, then $\{H_n\}_n$ is called a Hopf algebra filtration over
$R$.  The completion $\hat H = \plim_n H/H_n$ is a complete Hopf
algebra over $\hat R$ with the structure morphisms
\begin{equation*}
  \hat\Delta \colon \hat H \rightarrow  \hat H \hat\otimes  \hat H,\quad
  \hat\epsilon \colon \hat H \rightarrow  \hat R,\quad
  \hat S\colon \hat H \rightarrow  \hat H
\end{equation*}
induced by those of $H$.  Here $\hat H \hat\otimes  \hat H$ is the completed
tensor product of two copies of $\hat H$ defined by
\begin{equation*}
  \begin{split}
    \hat H \hat\otimes  \hat H
    &=\plim_{k,l} (\hat H\otimes _{\hat R}\hat H)/
  (\bar H_k\otimes _{\hat R}\hat H+\hat H\otimes _{\hat R}\bar H_l)\\
  &\simeq\plim_{k,l} (H\otimes _RH)/(H_k\otimes _RH+H\otimes _RH_l),
  \end{split}
\end{equation*}
where $\bar H_k$ is the closure of $H_k$ in $\hat H$.

Let $H$ be a Hopf algebra over a commutative ring $R$.  Let $I$ be an
ideal in $R$.  A two-sided ideal $J$ in $H$ is called a {\em Hopf
  ideal with respect to $I$} if
\begin{equation*}
  I\subset J,\quad \Delta (J)\subset J\otimes _RH+H\otimes _RJ,\quad
  S(J)\subset J,\quad \epsilon (J)\subset I.
\end{equation*}
In this case, the $J$-adic filtration $\{J^n\}_n$ of $H$ is a Hopf
algebra filtration with respect to $\{I^n\}_n$, and hence the
completion $\plim_n H/J^n$ is a complete Hopf algebra over $\hat
R=\plim_n R/I^n$.  Note that a Hopf ideal with respect to $(0)$ is
just a Hopf ideal in the usual sense.

The filtration \eqref{eq:F1} of $R$ will be called {\em fine} if
\begin{equation*}
  R_{m+n}\subset R_mR_n\quad (m,n\ge 0).
\end{equation*}
(Note the direction of the inclusion.)
The Hopf algebra filtration \eqref{eq:F2} will be called {\em fine} if
\begin{equation}
  \label{eq:F12} H_{m+n}\subset H_mH_n\quad (m,n\ge 0).
\end{equation}

Suppose that $\{R_n\}_n$ is a fine filtration of a commutative ring
$R$ with unit, and $\{H_n\}_n$ is a fine Hopf algebra filtration of a
Hopf algebra $H$ with respect to $\{R_n\}_n$.  Then the two-sided
ideal $H_1$ in $H$ is a Hopf ideal with respect to $R_1$.  Moreover,
for each $n\ge 0$,
\begin{equation*}
  R_n\subset R_1^n,\quad H_n\subset H_1^n.
\end{equation*}
Hence $\id_R$ induces a continuous ring homomorphism
\begin{equation*}
  \plim_n R/R_n \rightarrow  \plim_n R/R_1^n,
\end{equation*}
and $\id_H$ induces a continuous $(\plim_n R/R_n)$-algebra homomorphism
\begin{equation*}
  \plim_n H/H_n \rightarrow  \plim_n H/H_1^n.
\end{equation*}

\section{The filtration $\{\modU _n\}_n$ of $\modU $ and the completion $\hU$}
\label{sec:topology--}

\subsection{The filtration $\{\modA _n\}_n$ of $\modA $ and the completion $\hA$}
\label{sec:cycl-compl-ha}
For $n\ge 0$, set
\begin{gather*}
  (q)_n = (1-q)(1-q^2)\cdots (1-q^n)=(-1)^nv^{\hf n(n+1)}\{n\}!,\\
  \modA _n=(q)_n\modA =\{n\}!\modA ,\quad (\modA _q)_n = (q)_n\modA _q=\modA _n\cap \modA _q.
\end{gather*}
Since $\{m+n\}!=\{m\}!\{n\}!\bb{m+n}m\in \modA _m\modA _n$ for $m,n\ge 0$,
the filtrations $\{\modA _n\}_n$ and $\{(\modA _q)_n\}_n$ are fine.
Set
\begin{equation*}
  \hA = \plim_n \modA /\modA _n,\quad
  \hA_q = \plim_n \modA _q/(\modA _q)_n.
\end{equation*}
Note that the definition of $\hA$ is the same as in the introduction.
The $(\modZ /2)$-grading $\modA =\modA _q\oplus v\modA _q$ of $\modA $ induces the
$(\modZ /2)$-grading $\hA=\hA_q\oplus v\hA_q$ of $\hA$.  The ring $\hA_q$
is studied in \cite{H:cyclotomic}, where it was denoted by $\modZ [q]^\modN $.

\subsection{The filtration $\{\modU ^0_n\}_n$ of $\modU ^0$ and the completion
  $\hU^0$}
\label{sec:topology-0-0}

For each $n\ge 0$, let $\modU ^0_n$ denote the two-sided ideal in $\modU ^0$
generated by the elements $\BB{H+m}n$ for $m\in \modZ $.  Since
$\gamma _j(\BB{H+m}n)=\BB{H+m+j}n$ for each $j\in \modZ $, the ideal $\modU ^0_n$ of
$\modU ^0$ is preserved by $\gamma _j$, i.e.,
\begin{equation}
  \label{eq:1}
\gamma _j(\modU ^0_n) = \modU ^0_n.
\end{equation}
Hence, for any homogeneous $\modA $-submodule $M$ of $\modU $, we have
$M\modU ^0_n = \modU ^0_n M$, which we will freely use in what follows.
Another consequence of \eqref{eq:1} is that
$\gamma _j\colon \modU ^0\rightarrow \modU ^0$ induces a continuous $\hA$-algebra automorphism
\[
\gamma _j\colon \hU^0\rightarrow \hU^0.
\]

Since $\modU ^0=\modA [K,K^{-1}]$ is Noetherian, each ideal $\modU ^0_n$ is
finitely generated.  In fact, the following holds.

\begin{proposition}
  \label{thm:28}
  For each $n\ge 0$, the ideal $\modU ^0_n$ in $\modU ^0$ is generated by the
  elements
  \[
  \BB ni\BB H{n-i},\quad i=0,1,\ldots,n.
  \]
  In particular, for each $n\ge 0$ the ideal $\modU ^0_n$ contains the
  element $\{n\}!$, and hence we have $\{n\}!\modU ^0\subset \modU ^0_n$.
\end{proposition}

\begin{proof}
  Write $c_{n,i}=\BB ni\BB H{n-i}$, and set
  $(\modU ^0_n)'=(c_{n,0},c_{n,1},\ldots,c_{n,n})$.  We will show
  $\modU ^0_n=(\modU ^0_n)'$.  One can prove by induction that
  \begin{equation}
    \label{eq:T9}
    \BB{H+m}n = \sum_{i=0}^n v^{(n-i)m}K^{-i}\bb mi c_{n,i}\quad
    (m\in \modZ ,\ n\ge 0).
  \end{equation}
    Hence  $\modU ^0_n\subset (\modU ^0_n)'$.  To
  prove the other inclusion $(\modU ^0_n)'\subset \modU ^0_n$, we will see that
  $c_{n,p}\in \modU ^0_n$ for $p=0,\ldots,n$ by induction on $p$, the case
  $p=0$ being trivial.  If $0<p\le n$, then by \eqref{eq:T9} we have
  $c_{n,p}\in \BB{H+p}n + (c_{n,0},\ldots,c_{n,p-1})$, which is
  contained in $\modU ^0_n$ by the induction hypothesis.
\end{proof}

\begin{proposition}
  \label{thm:4}
  The family $\{\modU ^0_n\}_n$ is a fine Hopf algebra filtration of
  $\modU ^0$ with respect to $\{\modA _n\}_n$.
\end{proposition}

\begin{proof}
  We will verify the conditions \eqref{eq:F3}--\eqref{eq:F12} for $H_n=\modU ^0_n$
  and $R_n=\modA _n$.
  Proposition \ref{thm:28} implies \eqref{eq:F3}.
  One can verify by induction that
  \begin{equation*}
    \Delta (\BB {H+m}n) =\sum_{i+j=n} v^{-im}\bb ni K^{-i}
    {\BB{H+m}j} \otimes  K^{j} \BB Hi,
  \end{equation*}
  which implies \eqref{eq:F4}.  Since $\epsilon (\BB{H+m}n)=\BB mn=\bb
  mn\{n\}!\in \modA _n$, we have \eqref{eq:F5}.  Since
  $S(\BB{H+m}n)=\BB{-H+m}n=(-1)^n\BB{H-m+n-1}n\in \modU ^0_n$, we
  have \eqref{eq:F6}.
  Since $\BB{H+m}{n+n'}=\BB{H+m}n\BB{H+m-n}{n'}$ for $m\in \modZ $,
  $n,n'\ge 0$, we have \eqref{eq:F12}.
\end{proof}

Let $\hU^0$ denote the completion
\begin{gather*}
  \hU^0 = \plim_n \modU ^0/\modU ^0_n,
\end{gather*}
which is a complete Hopf $\hA$-algebra.

\subsection{A double filtration of $\modU ^0$ cofinal to $\{\modU ^0_n\}_n$}
\label{sec:double-filtration-0}

For $k,l\ge 0$, set
\begin{equation*}
  \modU ^0_{k,l} = (\{k\}!, \BB Hl)\subset \modU ^0.
\end{equation*}
The family $\{\modU ^0_{k,l}\}_{k,l\ge 0}$ forms a decreasing double
filtration of $\modU $, i.e., $\modU ^0_{k',l'}\subset \modU ^0_{k,l}$ if $0\le k\le k'$ and
$0\le l\le l'$.

\begin{proposition}
  \label{thm:20}
  The double filtration $\{\modU ^0_{k,l}\}_{k,l\ge 0}$ is cofinal to
  $\{\modU ^0_n\}_n$.
\end{proposition}

\begin{proof}
  By Proposition \ref{thm:28}, we have $\modU ^0_{n,n}=(\{n\}!,\BB
  Hn)\subset \modU ^0_n$ for $n\ge 0$, and $\modU ^0_{2n-1}\subset (\{n\}!,{\BB
    Hn})=\modU ^0_{n,n}$ for $n\ge 1$.  Therefore $\{\modU ^0_n\}_n$ is cofinal
  to $\{\modU ^0_{n,n}\}_n$, and hence to $\{\modU ^0_{k,l}\}_{k,l\ge 0}$.
\end{proof}

For $l\ge 0$, set
\[
{\BBB Hl} = (K^2-1)(K^2-q)\cdots(K^2-q^{l-1}) = v^{l(l-1)/2}K^l{\BB Hl}
\in \modA _q[K^2,K^{-2}].
\]

\begin{proposition}
  \label{thm:T2}
  There are natural isomorphisms
  \begin{equation}
    \label{eq:54}
    \hU^0
    \simeq \plim_l(\hA[K]/({\BBB Hl}))
    \simeq \plim_l(\hA[K,K^{-1}]/({\BB Hl})).
  \end{equation}
\end{proposition}

\begin{proof}
  By Proposition \ref{thm:20},
  \[
  \hU^0
  \simeq\plim_{k,l} \modU ^0/\modU ^0_{k,l}
  \simeq\plim_l(\plim_k \modA [K,K^{-1}]/(\{k\}!,\BBB Hl)).
  \]
  Note that $\BBB Hl\in \modA _q[K^2]\subset \modA [K]$ is a monic polynomial of
  degree $2l$ in $K$ with coefficients in $\modA $, and the degree $0$
  coefficient of $\BBB Hl$ is a unit in $\modA $.  Hence
  \[
  \begin{split}
    &\plim_k\modA [K,K^{-1}]/(\{k\}!,\BBB Hl)
    \simeq\plim_k\modA [K]/(\{k\}!,\BBB Hl)\\
    &\quad \simeq \plim_k((\modA /(\{k\}!))[K]/(\BBB Hl))
    \simeq \hA[K]/(\BBB Hl).
  \end{split}
  \]
  Therefore $ \hU^0\simeq\plim_l\hA[K]/(\BBB Hl)$.  The second
  isomorphism in \eqref{eq:54} is obvious from the above argument.
\end{proof}

Proposition \ref{thm:T2} implies the following.

\begin{corollary}
  \label{thm:6}
  Each element $t\in \hU^0$ is uniquely expressed as an infinite sum
  $t = \sum_{n\ge 0} t_n {\BBB Hl}$,
  where $t_n\in \hA+\hA K$ for $n\ge 0$.
\end{corollary}

\subsection{The filtration $\{\modU ^{\ge 0}_n\}_n$ of $\modU ^{\ge 0}$ and the
  completion $\hU^{\ge 0}$}
\label{sec:filtration-0nn-0}

For $n\ge 0$, let $\modU ^+_n$ denote the principal ideal in $\modU ^+_n=\modA [e]$
generated by the element $e^n$, and set
\begin{equation*}
  \modU ^{\ge 0}_n = \sum_{i+j=n} \modU ^0_i\modU ^+_j
  = \sum_{i+j=n} \modU ^+_j\modU ^0_i\subset \modU ^{\ge 0}.
\end{equation*}
Using $\modU ^{\ge 0}\modU ^0_i=\modU ^0_i\modU ^{\ge 0}$, we easily see that $\modU ^{\ge 0}$
is a two-sided ideal in $\modU ^{\ge 0}$.  Obviously, the family
$\{\modU ^{\ge 0}_n\}_n$ is a decreasing filtration.

\begin{proposition}
  \label{thm:1}
  The family $\{\modU ^{\ge 0}_n\}_n$ is a fine Hopf algebra filtration
  of $\modU ^{\ge 0}$ with respect to $\{\modA _n\}_n$.
\end{proposition}

\begin{proof}
  We will verify the conditions \eqref{eq:F3}--\eqref{eq:F12} for
  $H_n=\modU ^{\ge 0}_n$ and $R_n=\modA _n$ using Proposition \ref{thm:4}.  We
  clearly have \eqref{eq:F3}.  The inclusion \eqref{eq:F4} follows
  from Proposition \ref{thm:4} and
  \begin{equation*}
    \Delta (\modU ^+_n)\subset \sum_{i+j=n}\modU ^0\modU ^+_i\otimes _\modA \modU ^+_j,
  \end{equation*}
  which follows from \eqref{eq:h13}.  The condition \eqref{eq:F5}
  follows from Proposition \ref{thm:4} and the fact that
  $\epsilon (\modU ^+_n)=0$ for $n>0$.  We show \eqref{eq:F6} as follows:
  \[
  S(\modU ^{\ge 0}_n)=\sum_{i+j=n}S(\modU ^0_i\modU ^+_j)
  \subset \sum_{i+j=n}\modU ^0\modU ^+_j\modU ^0_i
  \subset \sum_{i+j=n}\modU ^0_i\modU ^+_j
  =\modU ^{\ge 0}_n.
  \]
  Here we used $S(\modU ^+_n)\subset \modU ^0\modU ^+_n$, which follows
  from \eqref{eq:40}.
  We see \eqref{eq:F12} as follows.  Suppose $m+n=i+j$, $i,j\ge 0$.  If
  $n\le j$, then
  $\modU ^0_i\modU ^+_j=\modU ^0_i\modU ^+_{j-n}\modU ^+_n\subset \modU ^{\ge 0}_m\modU ^{\ge 0}_n$,
  and otherwise
  $\modU ^0_i\modU ^+_j=\modU ^0_m\modU ^0_{i-m}\modU ^+_j\subset \modU ^{\ge 0}_m\modU ^{\ge 0}_n$.
  Hence \eqref{eq:F12} follows.
\end{proof}

Let $\hU^{\ge 0}$ denote the completion
\begin{equation*}
  \hU^{\ge 0} = \plim_n \modU ^{\ge 0}/\modU ^{\ge 0}_n,
\end{equation*}
which is a complete Hopf $\hA$-algebra.

Clearly, the filtration $\{\modU ^{\ge 0}_n\}_n$ is cofinal to the double
filtration $\{\modU ^0_k \modU ^++\modU ^0 \modU ^+_l\}_{k,l\ge 0}$.
We have
\[
  \modU ^{\ge 0}/(\modU ^0_k \modU ^++\modU ^0 \modU ^+_l)\simeq
  (\modU ^0\otimes _\modA \modU ^+)/(\modU ^0_k \modU ^++\modU ^0 \modU ^+_l)\simeq
  (\modU ^0/\modU ^0_k)\otimes _\modA (\modU ^+/\modU ^+_l).
\]
Therefore
\begin{equation*}
  \begin{split}
    \hU^{\ge 0}
    &\simeq \plim_{k,l} \modU ^{\ge 0}/(\modU ^0_k \modU ^++\modU ^0\modU ^+_l)
    \simeq \plim_{k,l} (\modU ^0/\modU ^0_k)\otimes _\modA (\modU ^+/\modU ^+_l)
  \end{split}
\end{equation*}
The right-hand side is isomorphic to $\plim_l \hU^0\otimes _\modA (\modU ^+/\modU ^+_l)$
since $\modU ^+/\modU ^+_l$ is a finitely generated free $\modA $-module. Hence
there is an isomorphism of complete $\hA$-modules
\begin{equation}
  \label{eq:132}
  \plim_k \hU^0\otimes _\modA (\modU ^+/\modU ^+_k) \simeqto \hU^{\ge 0},
\end{equation}
induced by $\modU ^0\otimes _\modA \modU ^+\simeqto\modU ^{\ge 0}$, $x\otimes y\mapsto xy$.

\begin{proposition}
  \label{thm:16}
  Each element $a\in \hU^{\ge 0}$ is uniquely expressed as an infinite sum
  $a = \sum_{n\ge 0} t_n e^n$,
  where $t_n \in \hU^0$ for $n\ge 0$.
\end{proposition}

\begin{proof}
  The result follows using \eqref{eq:132} since each element of
  $\plim_k \hU^0\otimes _\modA (\modU ^+/\modU ^+_k)$ is uniquely expressed as
  $a=\sum_{n\ge 0}t_n\otimes e^n$ with $t_n \in \hU^0$ for $n\ge 0$.
\end{proof}

\subsection{The filtration $\{\modU _n\}_n$ of $\modU $ and the completion $\hU$}
\label{sec:topology}

We have defined in the introduction the two-sided ideal $\modU _n$ in $\modU $
for $n\ge 0$.  We have
\begin{equation*}
  \modU _n = \modU \modU ^{\ge 0}_n\modU =\modU ^-\modU ^{\ge 0}_n\modU ^-.
\end{equation*}
Here the second identity follows from $\modU =\modU ^-\modU ^{\ge 0}=\modU ^{\ge 0}\modU ^-$.
  The family $\{\modU _n\}_n$ is a decreasing filtration of
two-sided ideals in $\modU $.  The following follows immediately from
Propositions \ref{thm:h1} and \ref{thm:1}.

\begin{proposition}
  \label{thm:5}
  The filtration $\{\modU _n\}_n$ is a fine Hopf algebra filtration
  of $\modU $ with respect to $\{\modA _n\}_n$.
\end{proposition}

The completion $\hU = \plim_n \modU /\modU _n$, already defined in
the introduction, is a complete Hopf $\hA$-algebra.

\subsection{The filtration $\{\modU '_n\}_n$}
\label{sec:filtration-nn0}

We will regard $\modU $ as a free right $\modU ^{\ge 0}$-module freely generated
by the elements $F^{(i)}$ for $i\ge 0$, the right action being
multiplication from the right.  For each $n\ge 0$, set
\[
  \modU '_n = \modU \modU ^{\ge 0}_n = \modU ^-\modU ^{\ge 0}_n,
\]
which is the left ideal in $\modU $ generated by $\modU ^{\ge 0}_n$, and is also
 a right $\modU ^{\ge 0}$-submodule of $\modU $.  The family
$\{\modU '_n\}_n$ is a decreasing filtration of left ideals in $\modU $.

\begin{proposition}
  \label{thm:18}
  The filtration $\{\modU '_n\}_n$ is cofinal to $\{\modU _n\}_n$.  In fact,
  \begin{equation}
    \label{eq:136}
    \modU _{2n-1}\subset \modU '_n\subset \modU _n\quad (n\ge 1).
  \end{equation}
\end{proposition}

\begin{proof}
  Obviously, $\modU '_n\subset \modU _n$ for $n\ge 0$.  By \eqref{eq:117},
  \begin{equation*}
    \modU ^+_m\modU ^- \subset \modU ^-\modU ^{\ge 0}_m
  \end{equation*}
  for all $m\ge 0$.  Hence for $i,j\ge 0$
  \[
  \modU ^-\modU ^0_i\modU ^+_j\modU ^-
  \subset \modU ^-\modU ^0_i\modU ^-\modU ^{\ge 0}_j
  =\modU ^-\modU ^0_i\modU ^{\ge 0}_j
  \subset \modU ^-\modU ^{\ge 0}_{\max(i,j)}
  =\modU '_{\max(i,j)}.
  \]
  Therefore
  \[
  \begin{split}
    &\modU _{2n-1}
    =\sum_{i+j=2n-1}\modU ^-\modU ^0_i\modU ^+_j\modU ^-
    \subset \sum_{i+j=2n-1}\modU '_{\max(i,j)}
    \subset \modU '_n.
  \end{split}
  \]
\end{proof}

By Proposition \ref{thm:18}, there is a natural isomorphism
\begin{equation*}
  \hU \simeq \plim_n\modU /\modU '_n = \plim_n\modU ^-\modU ^{\ge 0}/\modU ^-\modU ^{\ge 0}_n.
\end{equation*}
By Theorem \ref{thm:3}, we have
\[
\hU\simeq\plim_n\modU ^-\modU ^{\ge 0}/\modU ^-\modU ^{\ge 0}_n
\simeq\plim_n(\modU ^-\otimes _\modA \modU ^{\ge 0})/(\modU ^-\otimes _\modA \modU ^{\ge 0}_n)
\simeq\plim_n\modU ^-\otimes _\modA (\modU ^{\ge 0}/\modU ^{\ge 0}_n).
\]
Since $\modU ^-$ is freely $\modA $-spanned by the elements $F^{(l)}$ for
$l\ge 0$, the completion $\plim_n \modU ^-\otimes (\modU ^{\ge 0}/\modU ^{\ge 0}_n)$ is the
topologically-free right $\hU^{\ge 0}$-module generated by the
$F^{(l)}\otimes 1$ for $l\ge 0$.  Hence we have the following.

\begin{proposition}
  \label{thm:11}
  The algebra $\hU$ is a topologically-free right $\hU^{\ge 0}$-module
  which is
  topologically-freely generated by the elements $F^{(n)}$ for $n\ge 0$.
  Consequently, each element $a\in \hU$ is uniquely expressed as an
  infinite sum
  $a = \sum_{n\ge 0} F^{(n)} a_n$,
  where the sequence $a_n\in \hU^{\ge 0}$ for $n\ge 0$ converges to $0$ in
  $\hU^{\ge 0}$.
\end{proposition}

\begin{corollary}
  \label{thm:9}
  Each element $a\in \hU$ is uniquely expressed as an infinite sum
  $a = \sum_{m,n\ge 0} F^{(n)}t_{m,n} e^m$,
  where the elements $t_{m,n}\in \hU^0$ ($m,n\ge 0$) are such that for
  each $m\ge 0$ the sequence $(t_{m,n})_{n\ge 0}$ converges to $0$ in
  $\hU^0$ as $n\rightarrow \infty $.
\end{corollary}

\begin{proof}
  We express $a$ as in Proposition \ref{thm:11}.  By
  Proposition \ref{thm:16}, we can express each $a_n$ as
  $a_n = \sum_{m\ge 0} t_{m,n} e^m$,
  where $t_{m,n}\in \hU$ for each $m\ge 0$.  The assertion follows since
  $a_n$ converges to $0$ in $\hU^{\ge 0}$ as $n\rightarrow \infty $ if and only if for
  any $m\ge 0$ we have $t_{m,n}\rightarrow 0$ in $\hU^0$ as $n\rightarrow \infty $.
\end{proof}

\section{The $\modU _1$-adic topology for $\modU $ and the completion $\dU $}
\label{sec:U1-adic-topology}

\subsection{The $\modA _1$-adic completion $\hA$ of $\modA $}
\label{sec:_1-adic-completion}

Note that $\modA _1=(v-v^{-1})=(v^2-1)\subset \modA $, and $(\modA _q)_1=(q-1)\subset \modA _q$.
The completion $\dA=\plim_n\modA /\modA _1^n$ defined in the introduction is
$(\modZ /2)$-graded: $\dA\simeq \dA_q\oplus v\dA_q$, where we set
\begin{gather*}
  \dA_q=\plim_n\modA _q/(\modA _q)_1^n
  = \plim_n\modZ [q,q^{-1}]/(q-1)^n=\modZ [[q-1]].
\end{gather*}
As mentioned in the introduction, $\id_\modA $ induces an $\modA $-algebra
homomorphism $i_{\hA}\colon \hA\rightarrow \dA$.  Since $v-v^{-1}\in h\modQ [[h]]$,
$\id_\modA $ also induces an $\modA $-algebra homomorphism
$i_{\dA}\colon \dA\rightarrow \modQ [[h]]$.
\begin{proposition}
  \label{thm:13}
  $i_{\hA}$ and $i_{\dA}$ are injective.
\end{proposition}

\begin{proof}
  By \cite[Theorem 4.1]{H:cyclotomic}, the homomorphism
  $i_{\hA}|_{\hA_q}\colon \hA_q \rightarrow  \dA_q$ induced by $\id_{\modA _q}$ is
  injective.  (This result is also obtained by P. Vogel.)  The
  injectivity of $i_{\hA}$ follows from the $(\modZ /2)$-gradings of $\hA$
  and $\dA$.

  The homomorphism $i_{\dA}$ factors naturally as $\dA \rightarrow  \modZ [[v-1]] \rightarrow 
  \modQ [[h]]$.  Here $\dA=\plim_n\modZ [v,v^{-1}]/(v^2-1)^n\rightarrow \modZ [[v-1]]$ is
  injective by \cite[Corollary 4.1]{H:cyclotomic}, and
  $\modZ [[v-1]]\subset \modQ [[h]]$ is standard.
\end{proof}

In what follows, we regard $\hA\subset \dA\subset \modQ [[h]]$.

\subsection{The $\modU _1$-adic completion $\dU$ of $\modU $}
\label{sec:U1-adic-topology-1}

In view of Section \ref{sec:fine-hopf-algebra},
Propositions \ref{thm:4}, \ref{thm:1}, and \ref{thm:5}, we have
the following.

\begin{proposition}
  \label{thm:17}
  For $*=0,\ge 0,\none$, $\modU ^*_1$ is a Hopf ideal in $\modU ^*$ with respect
  to $\modA _1$.  The $\modU ^*_1$-adic filtration $\{(\modU ^*_1)^n\}_n$ is a
  fine Hopf algebra filtration with respect to the $\modA _1$-adic
  filtration $\{(\modA _1)^n\}_n$ of $\modA $.
\end{proposition}

Set
\begin{equation*}
  \dU^0 = \plim_n \modU ^0/(\modU ^0_1)^n,\quad
  \dU^{\ge 0} = \plim_n \modU ^{\ge 0}/(\modU ^{\ge 0}_1)^n,\quad
  \dU = \plim_n \modU /(\modU _1)^n,
\end{equation*}
which are complete Hopf $\hA$-algebras.
In view of section \ref{sec:fine-hopf-algebra}, the identity maps of
$\modU ^0$, $\modU ^{\ge 0}$, and $\modU $ induce
$\hA$-algebra homomorphisms
\begin{equation*}
  i_{\hU^0}\colon {\hU^0} \rightarrow  \dU^0,\quad
  i_{\hU^{\ge 0}}\colon {\hU^{\ge 0}} \rightarrow  \dU^{\ge 0},\quad
  i_{\hU}\colon {\hU} \rightarrow  \dU,
\end{equation*}
respectively.

\begin{proposition}
  \label{thm:25}
  For $n\ge 0$,
  \begin{gather}
    \label{eq:8} (\modU ^0_1)^n = \sum_{i=0}^n(\{1\}^i\{H\}^{n-i}),\\
    \label{eq:4} (\modU ^{\ge 0}_1)^n = \sum_{i+j=n}(\modU ^0_1)^i\modU ^+_j,\\
    \label{eq:3} (\modU _1)^n = \modU ^-(\modU ^{\ge 0}_1)^n,\\
    \label{eq:5} (\modU ^0_1)^n = (\modU _1)^n \cap  \modU ^0,\quad
    (\modU ^{\ge 0}_1)^n = (\modU _1)^n \cap  \modU ^{\ge 0}.
  \end{gather}
\end{proposition}

\begin{proof}
  In view of Proposition \ref{thm:28}, we have $\modU ^0_1=(\{1\},\{H\})$
  and hence \eqref{eq:8}.
  Since $\{H\}e^n=e^n\{H+2n\}=e^n(v^{2n}\{H\}+\{2n\}K^{-1})$, it
  follows that $\modU ^0_1e^n=e^n\modU ^0$, and hence
  $\modU ^0_1\modU ^+_n=\modU ^+_n\modU ^0_1$.  Using this identity, we
  have \eqref{eq:4}.
  By \eqref{eq:136}, $\modU _1=\modU '_1=\modU ^-\modU ^{\ge 0}_1$.  We also have
  $\modU _1=\modU ^{\ge 0}_1\modU ^-$, and hence $\modU ^{\ge 0}_1\modU ^-=\modU ^-\modU ^{\ge 0}_1$.
  Using this identity, we can easily verify \eqref{eq:3} using
  induction.
  \eqref{eq:5} follows from \eqref{eq:4}, \eqref{eq:3} and the
  triangular decomposition.
\end{proof}

In view of \eqref{eq:5}, the natural homomorphisms
$\dU^0\rightarrow \dU^{\ge 0}\rightarrow \dU$ are injective.  We will regard
$\dU^0\subset \dU^{\ge 0}\subset \dU$ in what follows.

The following follows from $\modU ^0_1=(v^2-1,K^2-1)\subset \modU ^0$.

\begin{proposition}
  \label{thm:36}
  We have
  \begin{equation*}
    \dU^0
    \simeq \plim_{k,l}\modU ^0/((v^2-1)^k,(K^2-1)^l)
    \simeq \modZ [[v^2-1, K^2-1]]\otimes _{\modZ [v^2,K^2]}\modZ [v,K].
  \end{equation*}
  Hence $\dU^0$ consists of the elements uniquely expressed as power
  series
  \begin{equation*}
    \sum_{i,j\ge 0}t_{i,j}(v^2-1)^i(K^2-1)^j,
  \end{equation*}
  where $t_{i,j}\in \modZ \oplus\modZ v\oplus\modZ K\oplus\modZ vK$ for $i,j\ge 0$.
\end{proposition}

Using arguments similar to those in Section
\ref{sec:filtration-0nn-0}, we obtain the following.
\begin{proposition}
  \label{thm:35}
  Each element of $\dU^{\ge 0}$ is uniquely expressed as an
  infinite sum
  $a = \sum_{n\ge 0} t_n e^n$,
  where $t_n\in \dU^0$ for $n\ge 0$.
\end{proposition}

  By an argument similar to Section \ref{sec:topology-0-0},
the completion $\dU$ is a topologically-free right $\dU^{\ge 0}$-module
topologically-freely generated by the elements $F^{(l)}$ for $l\ge 0$.
Hence we have the following.
\begin{proposition}
  \label{thm:21}
  Each element $a\in \dU$ is uniquely expressed as an infinite sum
  $a = \sum_{m,n\ge 0} F^{(n)}t_{m,n} e^m$,
  where the elements $t_{m,n}\in \dU^0$ ($m,n\ge 0$) are such that for
  each $m\ge 0$ the $t_{m,n}$ converges to $0$ in $\dU^0$ as $n\rightarrow \infty $.
\end{proposition}

\subsection{Embedding of $\modU $ into $\hU$ and $\dU$}
\label{sec:embedding--into}

\begin{proposition}
  \label{thm:37}
  We have
  $\bigcap_{n\ge 0}\modU _n =\bigcap_{n\ge 0}(\modU _1)^n = 0$.
  Consequently, the natural homomorphisms $\modU \rightarrow \hU$, $\modU \rightarrow \dU$ are
  injective.
\end{proposition}

\begin{proof}
  By $v^2-1, K^2-1, e\in hU_h$, we have
\begin{equation}
  \label{eq:121}
  \modU _1\subset hU_h.
\end{equation}
Using \eqref{eq:121} and the well-known fact that
$\bigcap_{n\ge 0}(hU_h)^n=0$, we conclude that
$\bigcap_{n\ge 0}(\modU _1)^n=0$, and hence the assertion follows.
\end{proof}

In what follows, we will regard $\modU \subset \hU$ and $\modU \subset \dU$ via the
injective homomorphisms in Proposition \ref{thm:37}.

\subsection{Embedding of $\dU $ into $U_h$}
\label{sec:embedding-du-}

By \eqref{eq:121}, there is an $\dA$-algebra homomorphism
\begin{equation*}
  i_{\dU}\colon \dU\rightarrow  U_h
\end{equation*}
induced by $\modU \subset U_h$.  We have the following.

\begin{proposition}
  \label{thm:22}
  The homomorphism $i_{\dU}$ is injective.
\end{proposition}

\begin{proof}
To prove injectivity of $i_{\dU}|_{\dU^0}\colon \dU^0\rightarrow U^0_h$, we factor it
as a composition of natural homomorphisms as follows.
\[
\begin{split}
 &\dU^0\simeq \plim_n\dA[K,K^{-1}]/(K^2-1)^n
\overset{i_1}\rightarrow \plim_n\modZ [[v-1]][K,K^{-1}]/(K^2-1)^n\\
&\quad \overset{i_2}\rightarrow  \modZ [[v-1,K-1]]
\overset{i_3}\rightarrow \modZ [[h,hH]]\overset{i_4}\rightarrow U_h^0.
\end{split}
\]
The map $i_1$ is injective since it is induced by $\dA\subset \modZ [[v-1]]$,
see the proof of Proposition \ref{thm:13}.  The map $i_2$ is injective
in view of \cite[Corollary 4.1]{H:cyclotomic}.  Injectivity of $i_3$
and $i_4$ is standard.  Therefore $i_{\dU}|_{\dU^0}$ is injective.

Injectivity of $i_{\dU}$ follows from injectivity of
$i_{\dU}|_{\dU^0}$ and Proposition \ref{thm:21}.
\end{proof}

In what follows, we regard $\dU $ as an $\dA$-subalgebra of $U_h$ via
the injective homomorphism $i_{\dU}$.

\subsection{Embedding of $\hU$ into $\dU $}
\label{sec:embedding-hb}

\begin{proposition}
  \label{thm:T7}
  For $*=0,\ge 0,\none$, the homomorphism $i_{\hU^*}$ is injective.
\end{proposition}

\begin{proof}
  First consider the case $*=0$.  Let $t\in \hU^0$ and suppose that
  $i_{\hU^0}(t) = 0$.  Express $t$ as in Proposition \ref{thm:6}.
  For $l\ge 1$, we can expand ${\BBB Hl}$ in $K^2-1$ as
  \[
  {\BBB Hl} = \prod_{i=0}^{l-1}((K^2-1) + (1-v^{2i}))
  = \sum_{j=1}^l b_{l,j} (v^2-1)^{l-j} (K^2-1)^j
  \]
  where $b_{l,j}\in \modA $ for $j=1,\ldots,l$.  We have $b_{l,l}=1$.  We have
  \[
  i_{\hU^0}(t)= t_0 +\sum_{l=1}^\infty  t_l
  \sum_{j=1}^l b_{l,j} (v^2-1)^{l-j} (K^2-1)^j
  = t_0 + \sum_{j=1}^\infty c_j(K^2-1)^j,
  \]
  where $c_j=\sum_{l=j}^\infty  t_l
  b_{l,j}(v^2-1)^{l-j}\in \dA+\dA K$.  Since $i_{\hU^0}(t)=0$,
  we have $t_0=0$, and, for each $j\ge 1$, we have $c_j=0$.  Since
  $b_{j,j}=1$, we can inductively verify that each $t_n$ $(n\ge 1)$ is
  divisible by $(v^2-1)^m$ for all $m\ge 0$.  We therefore have
  $t_1=t_2=\cdots=0$, consequently $t=0$.  Therefore $i_{\hU^0}$ is
  injective.

  The other two cases $*=\ge 0$ and $*=\none$ follow from the case
  $*=0$, Proposition \ref{thm:16}, Corollary \ref{thm:9},
  Proposition \ref{thm:35}, and Proposition \ref{thm:21}.
\end{proof}

In what follows, for $*=0,\ge 0,\none$, we regard $\hU^*$ as an
$\hA$-subalgebra of $\dU^*$ via $i_{\hU^*}$, and hence as an
$\hA$-subalgebra of $U_h^*$.

Since $\dU^0\subset U_h^0$ are integral domains, we have the following.

\begin{corollary}
  \label{thm:T9}
  The ring $\hU^0$ is an integral domain.
\end{corollary}

Since $\modU ^0_1$ is preserved by $\gamma _j$ for each $j\in \modZ $, there is an
induced automorphism $\gamma _j\colon  \dU^0\rightarrow \dU^0$.  In view of the inclusion
$\hU^0\subset \dU^0\subset U_h^0$, the automorphism $\gamma _j$ of $\hU^0$ and that of
$\dU^0$ are restrictions of $\gamma _j\colon U_h^0\rightarrow U_h^0$.

\section{The filtration $\{\modU ^e_n\}_n$, the completion $\cU$, and
  the subalgebra $\tU$ of $\hU$}
\label{sec:topology-c-}

\subsection{The filtration $\{\modU ^e_n\}_n$ of $\modU $ and the completion
  $\cU$ of $\modU $}
\label{sec:filtration-enn-}
Recall from the introduction that, for $n\ge 0$, $\modU ^e_n$ denotes the
two-sided ideal in $\modU $ generated by the element $e^n$.  We have
\begin{equation}
  \label{eq:53}
  \modU ^e_n = \modU e^n\modU  = \modU \modU ^+_n\modU  = \modU ^-\modU ^0\modU ^+_n\modU ^-=\modU ^0\modU ^-\modU ^+_n\modU ^-,
\end{equation}
Clearly, $\{\modU ^e_n\}_n$ is a decreasing filtration.

\begin{proposition}
  \label{thm:8}
  The family $\{\modU ^e_n\}_n$ is a fine Hopf algebra filtration of
  $\modU $ over $\modA $.
\end{proposition}

\begin{proof}
  The assertion follows from \eqref{eq:h13} and \eqref{eq:40}.
\end{proof}

The completion
$\cU = \plim_n \modU /\modU ^e_n$, defined in the introduction,
is a complete Hopf $\modA $-algebra.

\subsection{The image $\tU$}
\label{sec:image-tu}

We have $\modU ^e_n\subset \modU _n$ for $n\ge 0$.  Hence $\id_\modU $ induces an
$\modA $-algebra homomorphism
\begin{equation*}
  i_{\cU}\colon \cU \rightarrow  \hU,
\end{equation*}

\begin{conjecture}
  \label{thm:14}
  $i_{\cU}$ is injective.
\end{conjecture}

Let $\tU$ denote the image $i_{\cU }(\cU )\subset \hU$, which is an
$\modA $-subalgebra of $\hU$.

\subsection{Hopf algebra structure of $\cU$}
\label{sec:hopf-algebra-struct}

If Conjecture \ref{thm:14} holds, then we can identify $\tU$ and
$\cU$, and hence $\tU$ has a complete Hopf algebra structure.  Even if
Conjecture \ref{thm:14} is not true, we can define a Hopf algebra
structure on $\tU$ in a suitable sense as follows.  For $n\ge 0$, let
$\tU^{\tilde\otimes n}$ denote the image of the natural map
\[
  \cU^{\check\otimes n}\rightarrow \hU^{\hat\otimes n}
\]
between the completed tensor products.  We have natural isomorphisms
$\tU^{\tilde\otimes 0}\simeq\modA $ and $\tU^{\tilde\otimes 1}\simeq\tU$.  The
complete Hopf algebra structure of $\cU$ induces continuous
$\modA $-module maps
\begin{equation*}
  \tilde\Delta \colon \tU \rightarrow  \tU^{\tilde\otimes 2},\quad
  \tilde\epsilon \colon \tU \rightarrow  \modA ,\quad
  \tilde S\colon \tU \rightarrow  \tU.
\end{equation*}
These maps satisfy the Hopf algebra axioms, such as coassociativity
\[
(\tilde\Delta \otimes \id)\tilde\Delta =(\id\otimes \tilde\Delta )\tilde\Delta \colon \tU\rightarrow \tU^{\tilde\otimes 3}.
\]

\section{A $(\modZ /2)$-grading of $U$ and the induced gradings}
\label{sec:2-gradings-uq}

Let $\modG _0U$ denote the $\modQ (v)$-subalgebra of $U$ generated by the
elements $K^2$, $K^{-2}$, $E$, and $FK$.  Set $\modG _1U=K\modG _0U$.
Then we have a direct sum decomposition
\begin{equation}
  \label{eq:55}
  U = \modG _0U \oplus \modG _1U,
\end{equation}
which gives $U$ a $(\modZ /2)$-graded $\modQ (v)$-algebra structure.  An
element of $U$ is called {\em $K$-homogeneous} if it is homogeneous
in the grading \eqref{eq:55}.   Moreover, an additive subgroup $M$ of $U$
is called {\em $K$-homogeneous} if $M$ is homogeneous in the
grading \eqref{eq:55}, i.e., $M$ is generated by the $K$-homogeneous
element.  For any $K$-homogeneous additive subgroup $M$ of $U$ and
$i=0,1$, set
\begin{equation*}
  \modG _iM = \modG _iU \cap M.
\end{equation*}

It is easy to see that $\modU $ is $K$-homogeneous.  Hence \eqref{eq:55}
induces a $(\modZ /2)$-graded $\modA $-algebra structure on $\modU $:
\begin{equation*}
  \modU  = \modG _0\modU  \oplus \modG _1\modU ,
\end{equation*}
and $\modG _0\modU $ is generated as an $\modA $-algebra by the elements $K^2$,
$K^{-2}$, $e$, and $F^{(i)}K^i$ for $i\ge 1$.

The two-sided ideals $\modU _n$, $(\modU _1)^n$, and $\modU ^e_n$ for $n\ge 0$ are
all $K$-homogeneous since the generators of them are $K$-homogeneous.
Therefore $\vU=\hU, \dU, \cU, \tU$ inherits from $\modU $ a
$(\modZ /2)$-graded algebra structure
\begin{equation*}
  \vU=\modG _0\vU \oplus \modG _1\vU.
\end{equation*}
The degree $0$ part $\modG _0\vU$ is naturally identified with a
completion of $\modG _0\modU $ as follows:
\begin{gather*}
  \modG _0\hU \simeq \plim_n \modG _0\modU / \modG _0\modU _n,\quad
  \modG _0\dU \simeq \plim_n \modG _0\modU / \modG _0(\modU _1)^n,\\
  \modG _0\cU \simeq \plim_n \modG _0\modU / \modG _0\modU ^e_n,\quad
  \modG _0\tU \simeq i_{\cU}(\modG _0\cU).
\end{gather*}

\section{Centers}
\label{sec:centers}

\subsection{The Harish-Chandra homomorphism for $U_h$}
\label{sec:harish-chandra-homom}
Let
\begin{equation*}
  \varphi\colon \Gamma _0U_h\rightarrow U_h^0
\end{equation*}
denote the {\em Harish-Chandra homomorphism} of $U_h$, which is
defined to be the continuous $\modQ [[h]]$-module homomorphism satisfying
\begin{equation*}
  \varphi(F^itE^i) = \delta _{i,0}t
\end{equation*}
for $i\ge 0$ and $t\in U_h^0$.
Set
\begin{equation*}
  w=\gamma _2S=S\gamma _{-2}\colon U_h^0\rightarrow U_h^0,
\end{equation*}
which is the (involutive) continuous $\modQ [[h]]$-algebra automorphism of
$U_h^0$ determined by $w(H) = -H-2$.  Set
\begin{equation*}
  (U_h^0)^w= \{t\in U_h^0\ver w(t)=t\}.
\end{equation*}
We will use the following.

\begin{theorem}[\cite{T}]
  \label{thm:19}
  The homomorphism $\varphi$ maps the center $Z(U_h)$ isomorphically
  onto $(U_h^0)^w$.
\end{theorem}

\subsection{Facts about the centers of $\modU $, $\tU$, $\hU$, and $\dU $}
\label{sec:centers--hu}

In the rest of this section we will determine the centers of the
subalgebras $\modU \subset \tU\subset \hU\subset \dU$ of $U_h$.  Let $\vU$ denote one of
these algebras.  We have $Z(\vU) = \vU\cap  Z(U_h)$ since if $z\in \vU$ is
central in $\vU$, then $z$ commutes with $K,F,e$ and hence also with
$H$ and $E$, and consequently $z$ is central in $U_h$.

If $\vU=\modU $, $\hU$, or $\dU $, then let $\vU^0$ denote $\modU ^0$,
$\hU^0$, or $\dU ^0$, respectively.  In view of
Proposition \ref{thm:3}, Corollary \ref{thm:9}, and
Proposition \ref{thm:21}, we have
\begin{equation*}
  \vU^0 = \varphi(\Gamma _0\vU).
\end{equation*}
If $\vU=\tU$, then set
\begin{equation*}
  \vU^0=\tU^0 = \varphi(\Gamma _0\tU)\subset \hU^0.
\end{equation*}
({\it Warning\/}: the notation $\tU^0$ might suggest that $\tU^0=\tU\cap \hU^0$,
but this is not the case.  Actually, $\tU^0$ is not contained in
$\tU$.)

For $\vU=\modU ,\tU,\hU,\dU$, set
\begin{equation*}
  (\vU^0)^w = \vU^0\cap (U_h)^w = \{t\in \vU^0\ver w(t)=t\}.
\end{equation*}
In view of Theorem \ref{thm:19}, $\varphi$ maps $Z(\vU)$ injectively
into $(\vU^0)^w$.

\subsection{Center of $\modU $}
\label{sec:center-}
It is well-known that the center $Z(U)$ is freely generated as an
$\modQ (v)$-algebra by the element
\begin{equation*}
  C=fe+vK+vK^{-1} =(v-v^{-1})Fe +vK+v^{-1}K^{-1}\in  Z(\bU)\subset Z(\modU ),
\end{equation*}
i.e., we have
\begin{equation*}
  Z(U) = \modQ (v)[C].
\end{equation*}

\begin{theorem}
  \label{thm:29}
  The center $Z(\modU )$ of $\modU $ is freely generated as an $\modA $-algebra by
  the element $C$, i.e., we have
  \begin{equation}
    \label{eq:62}
    Z(\modU ) = \modA [C].
  \end{equation}
  The Harish-Chandra homomorphism $\varphi$ maps $Z(\modU )$
  isomorphically onto $(\modU ^0)^w$.  We also have
  \begin{equation*}
    Z(\modU )=Z(\bU).
  \end{equation*}
\end{theorem}

\begin{proof}
  Since $w(K^{\pm 1})=v^{\mp2}K^{\mp1}\in \modU $, the automorphism $w$ of
  $U_h^0$ restricts to an automorphism $w\colon \modU ^0\rightarrow \modU ^0$.  One can
  verify that the subalgebra $(\modU ^0)^w$ is freely generated as an
  $\modA $-algebra by the element $\varphi(C)=vK+vK^{-1}$.  By
  Section \ref{sec:centers--hu}, $Z(\modU )$ is mapped injectively into
  $(\modU ^0)^w$, and hence we have $Z(\modU )\subset \modA [C]$.  Since $C\in Z(\modU )$, we
  have $\modA [C]\subset Z(\modU )$.  Therefore we have \eqref{eq:62} and
  $\varphi(Z(U))=(\modU ^0)^w$.  The last assertion follows
  from \eqref{eq:62} and $C\in Z(\bU)$.
\end{proof}

Since $C\in \modG _1(Z(\modU ))$, it follows that the center $Z(\modU )$ is
$K$-homogeneous and
\begin{equation*}
  \modG _0Z(\modU ) = \modA [C^2],\quad
  \modG _1Z(\modU ) = C\modG _0Z(\modU )= C\modA [C^2].
\end{equation*}

\subsection{A basis of $Z(\modU )$}
\label{sec:basis-z}

As in the introduction, we set for $n\ge 0$
\begin{equation*}
  \sigma _n =\prod_{i=1}^n (C^2 - (v^i+v^{-i})^2)\in Z(\bU)\subset Z(\modU ).
\end{equation*}
Note that the $\sigma _n$ for $n\ge 0$ form a basis of $\modG _0Z(\modU )$, and hence
the $C\sigma _n$ for $n\ge 0$ form a basis of $\modG _1Z(\modU )$.

\begin{lemma}
  \label{lem:1}
  For $n\ge 0$,
  \begin{equation*}
    Z(\modU )\cap UE^nU \subset  \sigma _nZ(\modU ).
  \end{equation*}
\end{lemma}

\begin{proof}
  For $i\ge 1$, there is an $i$-dimensional irreducible left
  $U$-module $V^\pm _i$ generated by a highest weight vector
  $u^{\pm }_i$ satisfying $Ku^{\pm }_i=\pm v^{i-1}u^{\pm }_i$.  It is known
  that $C$ acts on $V^{\pm }_i$ by the scalar $\pm (v^i+v^{-i})$.

  Suppose that $z\in Z(\modU )\cap UE^nU$.  Write $z$ as a polynomial
  $z=g(C)\in \modA [C]$.  Since $z\in UE^nU$, $z$ acts as $0$ on the
  $i$-dimensional irreducible representations $V^+_i$ and $V^-_i$ with
  $i\le n$.  Since $z$ acts on $V^{\pm }_i$ by the scalar
  $g(\pm (v^i+v^{-i}))$, we have $g(\pm (v^i+v^{-i}))=0$ for $i=1,\ldots,n$.
  Hence the polynomial $z=g(C)$ is divisible by
  \[
  \prod_{i=1}^n(C-(v^i+v^{-i}))(C+(v^i+v^{-i}))
  =\prod_{i=1}^n(C^2-(v^i+v^{-i})^2)=\sigma _n.
  \]
\end{proof}

The following is proved in Section \ref{sec:proof-nen}.
\begin{proposition}
  \label{thm:23}
  For each $n\ge 0$, we have $\sigma _n\in \modU ^e_n$.
\end{proposition}

Assuming Proposition \ref{thm:23}, we have the following.

\begin{theorem}
  \label{thm:31}
  For each $n\ge 0$,
  \begin{equation*}
  \sigma _nZ(\modU ) = Z(\modU )\cap \modU ^e_n = Z(\modU )\cap UE^nU.
  \end{equation*}
\end{theorem}

\begin{proof}
  Lemma \ref{lem:1} and Proposition \ref{thm:23} implies that
  \[
  Z(\modU )\cap UE^nU\subset \sigma _nZ(\modU )\subset Z(\modU )\cap \modU ^e_n\subset   Z(\modU )\cap UE^nU,
  \]
  where the last inclusion is obvious.  Hence the assertion follows.
\end{proof}

\subsection{Center of $\hU$}
\label{sec:center-hu-2}

For $n\ge 0$, we have $w(\modU ^0_n)=\gamma _2S(\modU ^0_n)=\gamma _2(\modU ^0_n)=\modU ^0_n$.
Hence $w$ restricts to an involutive, continuous $\hA$-algebra
automorphism of $\hU^0$.  For $n\ge 0$, set
\begin{equation*}
    \bar\sigma _n
    = \varphi(\sigma _n)
    = \prod_{i=1}^n(\varphi(C)^2-(v^i+v^{-i})^2)
    = {\BB Hn}{\BB {H+1+n}n}.
\end{equation*}

\begin{lemma}
  \label{lem:2}
  We have $\hU^0\simeq \plim_n \hA[K,K^{-1}]/(\bar\sigma _n)$.
\end{lemma}

\begin{proof}
  Set $(\hU^0)'=\plim_n \hA[K,K^{-1}]/(\bar\sigma _n)$.  Since
  $(\bar\sigma _n)\subset (\BB Hn)$, $\id_{\hA[K,K^{-1}]}$ induces a continuous
  $\hA$-algebra homomorphism $f_1\colon (\hU^0)'\rightarrow \hU^0$.  We will show
  that there is a natural homomorphism $f_2\colon \hU^0\rightarrow (\hU^0)'$.  By an
  argument similar to the proof of Proposition \ref{thm:T2}, we see
  that $(\hU^0)'\simeq\plim_{m,n}\modA [K,K^{-1}]/(\{m\}!,\bar\sigma _n)$.
  Hence it suffices to prove that, for each $m,n\ge 0$, there is $n'\ge 0$
  such that $\BB H{n'}\in (\{m\}!,\bar \sigma _n)$.  By
  Proposition \ref{thm:28}, $\modU ^0_{4n}\subset (\BB{H}{2n+1},\{2n\}!)$.
  Applying the homomorphism $\gamma _{n+1}$ to both sides, we obtain
  $\modU ^0_{4n}\subset (\BB{H+n+1}{2n+1},\{2n\}!)\subset (\bar\sigma _n,\{2n\}!)$.  Hence,
  by Proposition \ref{thm:28}, we have $\BB{H}{n'}\in (\{m\}!,\bar\sigma _n)$
  where $n'=\max(4n,2m)$.  Hence $f_2$ exists.  Clearly, $f_2$ is
  inverse to $f_1$.
\end{proof}

\begin{remark}
  \label{thm:10}
  By the proof of Lemma \ref{lem:2}, we have
  \begin{equation*}
    \hU^0\simeq \plim_n\hA[K,K^{-1}]/(\BB{H+n+1}{2n+1}).
  \end{equation*}
  Therefore
  \begin{equation*}
    \begin{split}
      \hU^0
      &\simeq \plim_{X\subset \modZ ,\ |X|<\infty } \hA[K,K^{-1}] /
      (\prod_{k\in X}\{H+k\}),
    \end{split}
  \end{equation*}
  where $X$ runs through all the finite subsets of $\modZ $.
\end{remark}

Since $w\colon \hU^0\rightarrow \hU^0$ preserves the $(\modZ /2)$-grading
$\hU^0=\modG _0\hU^0\oplus\modG _1\hU^0$, it follows that $(\hU^0)^w$ is
$K$-homogeneous and hence has a $(\modZ /2)$-graded algebra structure:
$(\hU^0)^w=\modG _0(\hU^0)^w\oplus\modG _1(\hU^0)^w$.

\begin{lemma}
  \label{prop:T5}
  Each element of $\modG _0(\hU^0)^w$ is uniquely expressed as an infinite
  sum
  $t=\sum_{n\ge 0} t_n\bar\sigma _n$, where $t_n\in \hA$ for $n\ge 0$.
  We have $\modG _1(\hU^0)^w=\varphi(C)\modG _0(\hU^0)^w$.
\end{lemma}

\begin{proof}
  By Lemma \ref{lem:2}, each element $t\in \modG _0\hU^0$ is uniquely
  expressed as $t=\sum_{n\ge 0} (t_n+t'_nK^2)\bar\sigma _n$, where
  $t_n,t'_n\in \hA$, $n\ge 0$.  Suppose that $t\in \modG _0(\hU^0)^w$.  It
  suffices to show that $t'_n=0$ for $n\ge 0$.  We may assume without
  loss of generality that $t_n=0$ for each $n\ge 0$, and hence $t=\sum_n
  t'_nK^2\bar\sigma _n$.  The proof is by induction on $n$.  Suppose that
  $t_0=t_1=\cdots=t_{n-1}=0$ for $n\ge 0$.  Then $t\equiv t'_
  nK^2\bar\sigma _n\pmod{(\bar\sigma _{n+1})}$.  Since $t=w(t)$, we have
  $t'_nK^2\bar\sigma _n\equiv
  t'_nv^{-4}K^{-2}\bar\sigma _n\pmod{(\bar\sigma _{n+1})}$, and hence $t'_n(K^2-
  v^{-4}K^{-2})\in (\{H-n\}\{H+n+2\})$.  Since $K^2-v^{-4}K^{-2}$ is not
  divisible by $\{H-n\}\{H+n+2\}$, we have $t'_n=0$.  This completes
  the proof of the first assertion.

  Suppose that $t\in \modG _1(\hU^0)^w$.  By Lemma \ref{lem:2}, $t$ is uniquely
  expressed as an infinite sum $t=\sum_nt_n\bar\sigma _n$, where
  $t_0,t_1,\ldots\in \hA K+\hA K^{-1}$.  Then
  $w(t)=\sum_nw(t_n)\bar\sigma _n$, and $w(t_n)\in \hA K+\hA K^{-1}$ for
  $n\ge 0$.  Since $w(t)=t$, we have $\sum_n(t_n-w(t_n))\bar\sigma _n=0$, and
  hence $t_n=w(t_n)$ for $n\ge 0$.  Therefore
  $t_n\in \hA\varphi(C)$ for all $n\ge 0$, and hence
  $t\in \varphi(C)\modG _0(\hU^0)^w$.  This completes the second assertion.
\end{proof}

\begin{theorem}
  \label{thm:33}
  Each element $z\in Z(\hU)$ is uniquely expressed as an infinite sum
  $z = \sum_{n\ge 0} z_n \sigma _n$,
  where $z_n\in \hA +\hA C$ for $n\ge 0$.  The Harish-Chandra homomorphism
  $\varphi$ maps $Z(\hU)$ isomorphically onto $(\hU^0)^w$.
  We have
  \begin{equation*}
    Z(\hU) \simeq \plim_{m,n}Z(\modU )/
    (\{m\}!,\sigma _n)\simeq\plim_n\hA[C]/(\sigma _n).
  \end{equation*}
  The center $Z(\hU)$ is $K$-homogeneous, and we have
  \begin{equation*}
    \modG _0Z(\hU)\simeq \plim_n\hA[C^2]/(\sigma _n),\quad
    \modG _1Z(\hU) = C\modG _0Z(\hU).
  \end{equation*}
\end{theorem}

\begin{proof}
  By Section \ref{sec:centers--hu}, $\varphi$ maps $Z(\hU)$
  injectively into $(\hU^0)^w$.  By Lemma \ref{prop:T5}, each element
  $t$ of $(\hU^0)^w$ is uniquely expressed as an infinite sum $t =
  \sum_{n\ge 0}(t_n+t'_n\varphi(C))\bar\sigma _n$, where $t_n,t'_n\in \hA$ for
  $n\ge 0$.  The element $z = \sum_{n\ge 0}(t_n+t'_n C)\sigma _n$ is contained
  in $Z(\hU)$ since we have $\sigma _n\in \modU _n$ by Proposition \ref{thm:23}.
  Obviously, $\varphi(z)=t$.  Hence the first and the second
  assertions follow.  The rest of the theorem easily follows.
\end{proof}

\subsection{Center of $\tU$}
\label{sec:center-tu}

\begin{lemma}
  \label{lem:B1}
  If $n\ge 0$, then $\varphi(\Gamma _0\modU ^e_n) = (\BB Hn)\subset \modU ^0$.
  Consequently, each element of $\tU^0$ is uniquely expressed as an
  infinite sum $t=\sum_{n\ge 0}t_n\BB Hn$, where $t_n\in \modA $ for $n\ge 0$.
\end{lemma}

\begin{proof}
  In view of \eqref{eq:53}, $\Gamma _0\modU ^e_n$ is $\modA $-spanned by the
  elements
  \begin{equation*}
    t F^{(i)} e^m F^{(k)} e^l\quad (i,k,l\ge 0, m\ge n, m+l=i+k, t\in \modU ^0).
  \end{equation*}
  Since $\varphi(t F^{(i)} e^m F^{(k)} e^l)=0$ if either $i>0$ or
  $l>0$, it follows that $\varphi(\Gamma _0\modU ^e_n)$ is $\modA $-spanned by the
  elements $t \varphi(e^mF^{(m)})$ for $m\ge n$ and $t\in \modU ^0$.
  By \eqref{eq:117}, we have $\varphi(e^mF^{(m)})=\BB Hm$.  Hence
 the first assertion follows.
\end{proof}

\begin{lemma}
  \label{lem:8}
  Each element of $\modG _0(\tU^0)^w$ is uniquely expressed as an infinite
  sum
  $t= \sum_{n\ge 0} t_n \bar\sigma _n$,
  where $t_n\in \modA $ for $n\ge 0$.  We have
  $\modG _1(\tU^0)^w=\varphi(C)\modG _0(\tU^0)^w$.
\end{lemma}

\begin{proof}
  For $j\in \modZ $, let $s_j\colon \hU\rightarrow \hA$ denote the continuous $\hA$-algebra
  homomorphism determined by $s_j(K)=v^j$.

  Let $t$ be an element of $\modG _0(\tU^0)^w$ (resp. $\modG _1(\tU^0)^w$).
  By Lemma \ref{prop:T5}, $t$ is uniquely expressed as $t=\sum_{n\ge 0}
  t_n\bar\sigma _n=\sum_{n\ge 0} t_n\BB Hn\BB{H+n+1}n$, where $t_n\in \hA$
  (resp. $t_n\in \hA\varphi(C)$) for $n\ge 0$.  We prove by induction on
  $n$ that $t_n\in \modA $ (resp.  $t_n\in \modA \varphi(C)$) for $n\ge 0$.  Let
  $n\ge 0$ and assume that $t_i\in \modA $ (resp.  $t_i\in \modA \varphi(C)$) if
  $0\le i<n$.  By Lemma \ref{lem:B1}, we have $s_n(t)\in \modA $, and hence
  $\modA \ni s_n(t)=\sum_{i=0}^n s_n(t_i) \BB ni\BB{n+i+1}i$.  By the
  inductive hypothesis, $s_n(t_n)\BB nn\BB{2n+1}n\in \modA $.  In the case
  $t\in \modG _0(\tU^0)^w$, we see using Lemma \ref{lem:4} below that
  $s_n(t_n)=t_n\in \modA $.  In the case $t\in \modG _1(\tU^0)^w$, we have
  $t_n=b_n\varphi(C)$ for some $b_n\in \hA$.  Hence
  $s_n(t_n)=b_n(v^{n+1}+v^{-n-1})\in \modA $.  Since $v^{n+1}+v^{-n-1}$ is
  equal up to multiplication of a power of $v$ to the product of some
  cyclotomic polynomials in $q=v^2$, it follows from Lemma \ref{lem:4}
  below that $b_n\in \modA $, and therefore $t_n\in \modA \varphi(C)$.
\end{proof}

For $i\ge 1$, let $\Phi _i(q)\in \modZ [q]\subset \modA _q$ denote the $i$th cyclotomic
polynomial in $q$.

\begin{lemma}
  \label{lem:4}
  Suppose that $t\in \hA$ and $t\Phi (q)\in \modA $, where
  $\Phi (q)=\prod_{i\in \modN }\Phi _i(q)^{\lambda (i)}\in \modZ [q]$ is the product of some
  powers of cyclotomic polynomials in $q=v^2$.  Then $t\in \modA $.
\end{lemma}

\begin{proof}
  If $t\in \hA_q$, then it follows from \cite[Proposition
  7.3]{H:cyclotomic} that $t\in \modA _q$.  If $t\in v\hA_q$, then the first
  case implies that $t\in v\modA _q$.  Then the general case immediately
  follows.
\end{proof}

\begin{theorem}
  \label{thm:34}
  Each element $z\in Z(\tU)$ is uniquely expressed as an infinite sum $z
  = \sum_{n\ge 0} z_n \sigma _n$, where $z_n\in \modA  +\modA  C$ for $n\ge 0$.  The
  Harish-Chandra homomorphism $\varphi$ maps $Z(\tU)$ isomorphically
  onto $(\tU^0)^w$.  We have
  \begin{equation*}
    Z(\tU) \simeq \plim_nZ(\modU )/(\sigma _n).
  \end{equation*}
  The center $Z(\tU)$ is $K$-homogeneous, and we have
  \begin{equation*}
    \modG _0Z(\tU)\simeq \plim_n\modA [C^2]/(\sigma _n),\quad
    \modG _1Z(\tU) = C\modG _0Z(\tU).
  \end{equation*}
\end{theorem}

\begin{proof}
  Let $z\in Z(\tU)$.  By Lemma \ref{lem:8}, we have
  $\varphi(z)=\sum_{n\ge 0}(a_n+a'_n\varphi(C))\bar\sigma _n$ with
  $a_n,a'_n\in \modA $ for $n\ge 0$.  Set $z'=\sum_{n\ge 0}(a_n+a'_nC)\sigma _n$,
  which is an element of $Z(\tU)$ by Proposition \ref{thm:23}, and we
  have $\varphi(z)=\varphi(z')$.  Since $\varphi$ is injective, we
  have $z=z'$.  Therefore  the first and the second assertions follow.
  The rest of the theorem follows from the first assertion.
\end{proof}

\subsection{Center of $\dU $}
\label{sec:center-1}
For completeness, we state the following results.

\begin{theorem}
  \label{thm:12}
  Each element $z\in Z(\dU )$ is uniquely expressed as a formal power series
  $z = \sum_{n\ge 0} z_n (C^2-[2]^2)^n$,
  where $z_n\in \dA +\dA C$ for $n\ge 0$.  (Note that $C^2-[2]^2=\sigma _1$.)
  The Harish-Chandra homomorphism $\varphi$ maps $Z(\dU)$
  isomorphically onto $(\dU^0)^w$.  We have
  \begin{gather*}
    Z(\dU )
    \simeq \plim_{m,n}Z(\modU )/(\{1\}^m,(C^2-[2]^2)^n)
    \simeq\plim_n\dA[C]/((C^2-[2]^2)^n),\\
    \modG _0Z(\dU)\simeq\dA[[C^2-[2]^2]],\quad
    \modG _1Z(\dU)= C\modG _0Z(\dU).
  \end{gather*}
\end{theorem}

The proof of Theorem \ref{thm:12} is similar to that of
Theorem \ref{thm:33}.  One also has a version of Theorem \ref{thm:12}
in which the occurrences of $(C^2-[2]^2)^n$ are replaced with $\sigma _n$
as follows.

\begin{theorem}
  \label{thm:15}
  Each element $z\in Z(\dU )$ is uniquely expressed as an infinite sum
  $z = \sum_{n\ge 0} z_n \sigma _n$, where $z_n\in \dA +\dA C$ for $n\ge 0$.  We
  have
  \begin{equation*}
    Z(\dU )
    \simeq \plim_{m,n}Z(\modU )/(\{1\}^m,\sigma _n)
    \simeq\plim_n\dA[C]/(\sigma _n).
  \end{equation*}
\end{theorem}

\begin{proof}
  Use Theorem \ref{thm:12} and the fact that the double filtrations
  $\{(\{1\}^m, \sigma _n)\}_{m,n}$ and $\{(\{1\}^m, (C^2-[2]^2)^n)\}_{m,n}$
  of $Z(\modU )$ are cofinal.
\end{proof}

\subsection{Multiplication in $Z(\modU )$, $Z(\hU)$, and $Z(\tU)$}
\label{sec:multiplication--hu}
Multiplication in the centers $Z(\modU )$, $Z(\hU)$, and $Z(\tU)$ is
described by giving a formula for the product $\sigma _m\sigma _n$ as a linear
combination of the $\sigma _j$.

\begin{proposition}
  \label{thm:26}
  If $m,n\ge 0$, then
  \begin{equation*}
    \begin{split}
    \sigma _m\sigma _n
    &= \sum_{i=0}^{\min(m,n)} \BB mi\BB ni\bb{m+n+1}i \sigma _{m+n-i}\\
    &= \sum_{j=\max(m,n)}^{m+n} \BB m{m+n-j}\BB n{m+n-j}\bb{m+n+1}{j+1} \sigma _{j}.
    \end{split}
  \end{equation*}
\end{proposition}

\begin{proof}
  The proof is a straightforward induction on $n$ using an identity
  \begin{equation*}
    \sigma _m\sigma _{n+1}
     = \sigma _{m+1}\sigma _n+\{m-n\}\{m+n+2\}\sigma _m\sigma _n.
  \end{equation*}
\end{proof}

\begin{remark}
  \label{thm:30}
  Suppose that $a = \sum_{m\ge 0}a_m\sigma _m$ and $b=\sum_{n\ge 0}b_n\sigma _n$ are
  elements of $\modG _0Z(\hU)$ (or $\modG _0Z(\tU)$, $\modG _0Z(\modU )$), then
  \begin{equation*}
    \begin{split}
      ab
      &= \sum_{j\ge 0}\Bigl(\sum_{0\le m,n\le j;\  m+n\ge j}  \BB{m}{m+n-j} \BB{n}{m+n-j}
      \bb{m+n+1}{j+1}a_mb_n\Bigr)\sigma _j.
    \end{split}
  \end{equation*}
\end{remark}

\section{Proof of Proposition \ref{thm:23}}
\label{sec:proof-nen}
In this section, we will prove
Proposition \ref{thm:23}.  In the first two subsections
\ref{sec:some-formulae},
and \ref{sec:adjoint-action}, we prepare necessary facts.  The proof
is given in \ref{sec:proof-prop-refthm:23}, in which we need a result
in \ref{sec:lemma-divisibility}.

\subsection{Formulas in $\modU ^0$}
\label{sec:some-formulae}

For $a\in \modZ H+\modZ $ and $r\ge 0$, we have the following well-known formula
\begin{equation}
  \label{eq:181}
  \{a\}_r = \sum_{j=0}^r (-1)^j v^{\hf(r-2j)(-r+1+2a)}\bb rj,
\end{equation}
which one can verify by induction on $r$.

For $r,s\ge 0$, $a,b,c\in \modZ H+\modZ $, set
\begin{equation*}
  \begin{split}
  \kappa (a, r; b, s; c)
  =&
  \sum_{i=0}^{r}\sum_{j=0}^{s} (-1)^{i+j}
  v^{\hf(r-2i)(-r+1+2a)+\hf(s-2j)(-s+1+2b) +\hf c(rs-4ji)}
  {\bb ri}{\bb sj}\\
  =&
  \sum_{j=0}^{s} (-1)^j
  v^{\hf(s-2j)(-s+1+2b+rc)}{\bb sj}{\BB {a+jc}r}.
  \end{split}
\end{equation*}
where the identity follows from \eqref{eq:181}.  We have the symmetry
\begin{equation}
  \label{eq:184}
  \kappa (a,r;b,s;c)=  \kappa (b,s;a,r;c).
\end{equation}
Note that $\kappa (a,r; b,s; c)\in v^{\hf crs}\modU ^0$.

\begin{lemma}
  \label{lem:6}
  If $0\le j\le m$, then
  \begin{equation}
    \label{eq:85}
    \sum_{i=j}^m(-1)^iv^{i(j-m)}{\bb ij}{\bb{m+1}{i+1}} = (-1)^jv^{m-j}.
  \end{equation}
\end{lemma}

\begin{proof}
  It is well-known that
  \begin{equation}
    \label{eq:150}
    \sum_{i=k}^l (-1)^iv^{i(k-l)}\bb{i-1}{k-1} \bb li = (-1)^k
  \end{equation}
  for $1\le k\le l$, which can also be derived from \eqref{eq:184} by
  setting $(a,r,b,s,c)=(-1,k-1,l-k,l,1)$.  The formula \eqref{eq:85}
  follows from \eqref{eq:150} by setting $l=m+1$, $k=j+1$, and
  shifting $i$ by $1$.
\end{proof}

\begin{lemma}
  \label{lem:7}
  If $n\in \modZ $ and $l\ge 0$, then
  \begin{equation}
    \label{eq:107}
    \sum_{j=0}^l(-1)^jv^{j(l-1-n)}{\bb lj}\BB{H-n}{l-j}\BB{H+j-l}j
    = (-1)^lv^{l^2-l}{\BB nl}K^{-l}.
  \end{equation}
\end{lemma}

\begin{proof}
  Let $\lambda _l$ and $\lambda '_l$ denote the left and right hand sides
  of \eqref{eq:107}, respectively.  It is straightforward to verify
  that the $\lambda _n$ and the $\lambda '_n$ satisfy the same recurrence relations:
  \[
  \begin{split}
  \lambda _{l+1}=\{H-n\}\gamma _{-1}(\lambda _l)-v^{2l-n}\{H-l\}\lambda _l,\\
  \lambda '_{l+1}=\{H-n\}\gamma _{-1}(\lambda '_l)-v^{2l-n}\{H-l\}\lambda '_l.
  \end{split}
  \]
  An induction proves the assertion.
\end{proof}

\subsection{Adjoint action}
\label{sec:adjoint-action}

Let $\trr\colon U\otimes U\rightarrow U$, $x\otimes y\mapsto x\trr y$, denote the (left)
adjoint action defined by
\begin{equation*}
  x\trr y = \sum x_{(1)} y S(x_{(2)})
\end{equation*}
for $x,y\in U$, where $\Delta (x)=\sum x_{(1)}\otimes x_{(2)}$.  We will regard
$U$ as a left $U$-module by the adjoint action.

For a homogeneous element $x$ in $U$ and $m\ge 0$, we have
\begin{gather}
  \label{eq:96}
  e^m\trr x = \sum_{i=0}^m(-1)^iv^{i(m-1+2|x|)}{\bb mi}
  e^{m-i}xe^i,\\
  \label{eq:100}
  F^{(m)}\trr x = \sum_{i=0}^m (-1)^iv^{-i(m-1)}
  F^{(m-i)}xF^{(i)}K^m,
\end{gather}
which can easily be derived using formulas in
Section \ref{sec:hopf-algebra-struct-2}.

For each $n\ge 0$, set
\begin{equation*}
  M_n=U\trr K^{-n}E^{(n)}=U\trr K^{-n}e^n \subset U,
\end{equation*}
which is an irreducible $(2n+1)$-dimensional left $U$-module, with a
highest weight vector $K^{-n}E^{(n)}$ of weight $v^{2n}$, i.e., we
have $E\trr K^{-n}E^{(n)}=0$ and $K \trr K^{-n}E^{(n)}=v^{2n}$.  We
have
\begin{equation*}
  F^{(2n)}\trr K^{-n}E^{(n)} = (-1)^nv^{-n(n-1)} F^{(n)},
\end{equation*}
and $F^{(k)}\trr K^{-n}E^{(n)}=0$ if $k>2n$.  The element $F^{(n)}$ is
a lowest weight vector, i.e., $F\trr F^{(n)}=0$.

\begin{lemma}
  \label{lem:3}
  If $n,j\ge 0$, then
  \begin{equation*}
    e^j\trr F^{(n)}
    =\sum_{k=0}^n v^{-\hf(j-k)(-j+1+2n)}{\bb jk}
    \beta _{n,j,k}
    F^{(n-k)}e^{j-k},
  \end{equation*}
  where
  \begin{equation*}
    \begin{split}
      \beta _{n,j,k}
      &=\kappa (n-k,j-k; H-j+n , k; 1)
      \\
      &=\sum_{i=0}^{j-k}
      (-1)^iv^{\hf(j-k-2i)(-j+1+2n)}{\bb {j-k}{i}}
      \BB{H-j+i+n}k\\
      &=\sum_{l=0}^{k} (-1)^l
      v^{\hf(k-2l)(1-j+2n-2k)}{\bb kl}\BB{n-k+l}{j-k}K^{k-2l}.
    \end{split}
  \end{equation*}
\end{lemma}

\begin{proof}
  The proof is by computation and we omit it.
\end{proof}

\subsection{Proof of Proposition \ref{thm:23}}
\label{sec:proof-prop-refthm:23}

For each $m\ge 0$, set
\begin{gather*}
  \xi _m  = \prod_{i=1}^m(C-v^{2i+1}-v^{-2i-1})\in Z(\bU),\\
  \xi '_m = \prod_{i=1}^m(C-v^{2i+1}K-v^{-2i-1}K^{-1})\in \Gamma _0\bU.
\end{gather*}
By induction on $m$, we can check that
  \begin{equation}
    \label{eq:83}
    \xi '_m = \{m\}!\sum_{i=0}^m (-1)^i {\BB{H+m+1}i} F^{(m-i)}e^{m-i}.
  \end{equation}
For $n\ge 0$, set
\begin{equation*}
  \wt\sigma _n
  = (\{2n\}!)^{-1}\xi _{2n}\trr F^{(n)}K^{-n}e^n
  = (\{2n\}!)^{-1}\xi '_{2n}\trr F^{(n)}K^{-n}e^n.
\end{equation*}
Here the identity holds since $K^{\pm 1}\trr x=x$ for any $x\in \Gamma _0\modU $.
In view of \eqref{eq:83}, we have $\wt\sigma _n\in \modU ^e_n$.

\begin{lemma}
  \label{lem:11}
  For each $n\ge 0$, we have $\wt\sigma _n\in \sigma _nZ(\modU )$.
\end{lemma}

\begin{proof}
  In view of Lemma \ref{lem:1} and $\wt\sigma _n\in \modU ^e_n$, it suffices to
  prove that $\wt\sigma _n$ is central.  The elements $F^{(n)}$ and
  $K^{-n}e^n$ are contained in the $(2n+1)$-dimensional, irreducible
  left $U$-submodule $M_n=U\trr K^{-n}E^{(n)}$ of the adjoint
  representation of $U$.  Hence the product $K^{-n}e^n$ is contained
  in $M_nM_n=\mu (M_n\otimes M_n)\subset U$, where $\mu \colon U\otimes U\rightarrow U$ is the
  multiplication.  Since $M_n$ is a $(2n+1)$-dimensional irreducible
  representation, it follows from the Clebsch-Gordan rule that
  $M_n\otimes M_n$ has a direct sum decomposition
  \begin{equation*}
    M_n\otimes M_n = W_0\oplus W_1\oplus\cdots\oplus W_{2n},
  \end{equation*}
  where, for each $i=0,\ldots,2n$, $W_i$ is the $(2n+1)$-dimensional
  irreducible $U$-module generated by a highest weight vector of
  weight $v^{2n}$.  Therefore
  \begin{equation*}
    M_nM_n = W'_0\oplus W'_1\oplus\cdots\oplus W'_{2n},
  \end{equation*}
  where $W'_i=\mu (W_i)$ is either isomorphic to $W_i$ or zero for each
  $i=0,\ldots,2n$.  (Actually, all the $W'_i$ are nonzero, but we will
  not need this fact.)  Since, for each $i=1,\ldots,m$, the element
  $C-v^{2i+1}-v^{-2i-1}$ acts as $0$ on $W_i$, the element $\xi _{2n}$
  acts by the adjoint action as $0$ on $W'_i$.  Therefore we have
  $\wt\sigma _n\in W_0'$, and hence $\wt\sigma _n$ is central.
\end{proof}

\begin{lemma}
  \label{lem:9}
  For $n\ge 0$, we have $\wt\sigma _n\in \modA \sigma _n$.
\end{lemma}

\begin{proof}
  By computation (see below), we obtain
  \begin{equation}
    \label{eq:177}
    \varphi(\wt\sigma _n)
    =(-1)^n v^{-n^2+n}\BB Hn
    \sum_{i=n}^{2n}(-1)^i \BB in \bb{2n+1}{i+1} \BB{H+n-i}n.
  \end{equation}
  It follows that
  \begin{equation*}
    \varphi(\wt\sigma _n)\in \Span_\modA \{K^{2n}, K^{2n-2},\ldots,K^{-2n}\}.
  \end{equation*}
  By Lemma \ref{lem:11}, we have
  $\wt\sigma _n\in \Span_\modA \{1,C^2,C^4,\ldots,C^{2n}\}$, and the assertion follows.

  Here we sketch the proof of \eqref{eq:177}.  It follows
  from \eqref{eq:83} and the fact that $t\trr x=\epsilon (t)x$ for $t\in \modU ^0$,
  $x\in \Gamma _0\modU $, we have
  \begin{equation}
    \label{eq:106}
    \begin{split}
      \wt\sigma _n
      =&\sum_{i=0}^{2n}(-1)^i\BB{2n+1}{2n-i}F^{(i)}e^i\trr F^{(n)}K^{-n}e^n.
    \end{split}
  \end{equation}
  If $i<n$, then $\varphi(F^{(i)}e^{i}\trr F^{(n)}K^{-n}e^n)=0$.
  Otherwise, we can prove using \eqref{eq:96} and \eqref{eq:100} that
  \[
  \begin{split}
    &\varphi(F^{(i)}e^{i}\trr F^{(n)}K^{-n}e^n)\\
    =&(-1)^i v^{-i(i-1)+2ni-2n^2}
    {\bb in}{\BB Hn}\BB{H+n-i}n K^{i-n}\\&\quad \quad
    \sum_{j=0}^{i-n}(-1)^j v^{j(i-1)-2nj}{\bb {i-n}j}
   \BB{H-n}{i-n-j}\BB{H+n+j-i}{j}.
  \end{split}
  \]
  By Lemma \ref{lem:7}, the sum in the right-hand side is equal to
  \[
  (-1)^{i-n}v^{(i-n)^2-(i-n)}\BB n{i-n}K^{-i+n}.
  \]
  Then it is easy to check \eqref{eq:177}.
\end{proof}

\begin{proposition}
  \label{thm:32}
  For $n\ge 0$, we have $\wt\sigma _n = v^{-n^2+n}\{n\}!\sigma _n$.
\end{proposition}

\begin{proof}
  For an element $t\in \modU ^0=\modA [K,K^{-1}]$ and $p\in \modZ $, let $c_p(t)\in \modA $
  denote the coefficient of $K^p$ in $t$.  Note that
  $c_n(\BB{H+r}n)=v^{\sum_{s=r-n+1}^{r}s}=v^{n(2r-n+1)/2}$ for
  $r\in \modZ $.  Hence
  \[
  c_{2n}(\varphi(\sigma _n))=c_{2n}({\BB Hn}\BB{H+n+1}n)
  =c_n({\BB Hn})c_n(\{H+n+1\})
  =v^{2n}.
  \]
  Using \eqref{eq:177}, we obtain
  \[
  c_{2n}(\varphi(\wt\sigma _n))
  =(-1)^n v^{-n^2+2n}\{n\}!\sum_{i=n}^{2n}
    (-1)^iv^{-ni}\bb in{\bb{2n+1}{i+1}}.
  \]
  By Lemma \ref{lem:6},
  \[
    c_{2n}(\varphi(\wt\sigma _n))
    =(-1)^n v^{-n^2+2n}\{n\}!(-1)^nv^n
     =v^{-n^2+3n}\{n\}!.
  \]
  Hence, by Lemma \ref{lem:9},
  \[
  \wt\sigma _n
  =\frac{c_{2n}(\varphi(\wt\sigma _n))}{c_{2n}(\varphi(\sigma _n))}\sigma _n
  = \frac{v^{-n^2+3n}\{n\}!}{v^{2n}}\sigma _n=v^{-n^2+n}\{n\}!\sigma _n.
  \]
\end{proof}

\begin{lemma}
  \label{lem:13}
  For $n\ge 0$,
  \begin{equation*}
    \wt\sigma _n
    =\sum_{j=0}^{2n}(-1)^j
    v^{j(-j-1+2n)+2n}\{2n-j\}!
    (e^j\trr F^{(n)})(F^{(j)}\trr K^{-n}e^n).
  \end{equation*}
\end{lemma}

\begin{proof}
  Since $e^iF^{(i)}-F^{(i)}e^i\in \modU \{H\}$ by \eqref{eq:117}, we have
  $F^{(i)}e^i\trr x= e^iF^{(i)}\trr x$ for $x\in \Gamma _0\modU $.  Hence it
  follows from \eqref{eq:117}, \eqref{eq:h13}, \eqref{eq:h5}, and
  $e\trr K^{-n}e^n=F^{(p)}\trr F^{(r)}=0$ ($p\ge 1$) that
  \[
  \begin{split}
    F^{(i)}e^i\trr F^{(n)}K^{-n}e^n
    =&e^iF^{(i)}\trr F^{(n)}K^{-n}e^n\\
    =&e^i \trr F^{(n)}(F^{(i)}\trr K^{-n}e^n)\\
    =&\sum_{j=0}^i v^{-j(i-j)}{\bb ij}
    (K^{i-j}e^j\trr F^{(n)})(e^{i-j}F^{(i)}\trr K^{-n}e^n).
  \end{split}
  \]
  Since
  \[
    K^{i-j}e^j\trr F^{(n)} = v^{2(i-j)(j-n)}e^j\trr F^{(n)}
  \]
  and
  \[
  e^{i-j}F^{(i)}\trr K^{-n}e^n
  =F^{(j)}\BB{H-j}{i-j}\trr K^{-n}e^n
  =F^{(j)}\BB{2n-j}{i-j}\trr K^{-n}e^n,
  \]
  we have
  \[
    F^{(i)}e^i\trr F^{(n)}K^{-n}e^n
    =\sum_{j=0}^i v^{(i-j)(j-2n)}{\bb ij}\BB{2n-j}{i-j}
    (e^j\trr F^{(n)})(F^{(j)}\trr K^{-n}e^n).
  \]
  By \eqref{eq:106},
  \[
  \begin{split}
    &\wt\sigma _n\\
    =&\sum_{i=0}^{2n}(-1)^i\BB{2n+1}{2n-i}
    \sum_{j=0}^i v^{(i-j)(j-2n)}{\bb ij}\BB{2n-j}{i-j}
    (e^j\trr F^{(n)})(F^{(j)}\trr K^{-n}e^n)\\
    =&\sum_{j=0}^{2n} v^{-j(j-2n)}\{2n-j\}!
    \Bigl(\sum_{i=j}^{2n}
    (-1)^i v^{i(j-2n)}{\bb ij}{\bb{2n+1}{i+1}}
     \Bigr)
    (e^j\trr F^{(n)})(F^{(j)}\trr K^{-n}e^n).
  \end{split}
  \]
  Here we used the identity
  $\BB{2n-j}{i-j}\BB{2n+1}{2n-i}=\{2n-j\}!{\bb{2n+1}{i+1}}$.  By
  Lemma \ref{lem:6}, the sum over $i$ is equal to $(-1)^jv^{2n-j}$.
  Hence the assertion follows.
\end{proof}

\begin{proof}[Proof of Proposition \ref{thm:23}]
  We have to show that $\wt\sigma _n\in \{n\}!\modU ^e_n$.  In view of
  Lemma \ref{lem:13} and $F^{(j)}\trr K^{-n}e^n\in \modU ^e_n$, it suffices
  to prove that
  \begin{equation*}
     \{2n-j\}!(e^j\trr F^{(n)})\in \{n\}!\modU 
  \end{equation*}
  for each $j=0,\ldots,2n$.  Note that the case $0\le j\le n$ is trivial.  We
  assume $n<j\le 2n$ in the following.  In view of Lemma \ref{lem:3}, it
  suffices to show that if $0\le l\le k\le n<j\le 2n$, then
  \begin{equation}
    \label{eq:183}
    \{2n-j\}!\bb jk \bb kl \BB{n-k+l}{j-k} \in \{n\}!\modA .
  \end{equation}
  Note that $\BB{n-k+l}{j-k}=0$ unless $l-j+n\ge 0$.  If $l-j+n\ge 0$,
  then \eqref{eq:183} is equal to
  \[
  \{n\}!\{n-k\}!\theta (j-n,k-l,l-j+n,n-k),
  \]
  which is contained in $\{n\}!\{n-k\}!\modA \subset \{n\}!\modA $ by
  Lemma \ref{lem:16} below.  This completes the proof.
\end{proof}

\subsection{An integrality lemma}
\label{sec:lemma-divisibility}

For $a,b,x,y\ge 0$, we set
\begin{equation*}
  \theta (a,b,x,y) = \psi (a,x,y){\bb{x+y}x}{\bb{b+x+y}b}{\bb{2a+b+x+y}a}\in \modQ (v)
\end{equation*}
where
\begin{equation*}
  \psi (a,x,y) = \frac{\{a+x+y\}!\{a\}!}{\{a+x\}!\{a+y\}!}\in \modQ (v).
\end{equation*}

\begin{lemma}
  \label{lem:16}
  If $a,b,x,y\ge 0$, then $\theta (a,b,x,y) \in \modA $.
\end{lemma}

\begin{proof}
  For each $n\ge 1$, let $\phi _n\in \modA $ denote the ``balanced cyclotomic
  polynomial''
  \[
  \phi _n = v^{-\deg\Phi _n(q)}\Phi _n(q),
  \]
  where $\Phi _n(q)\in \modZ [q]$ is the $n$th cyclotomic polynomial.  We have
  \begin{equation*}
    \{n\}! = \prod_{d\ge 1} \phi _d^{\fd{n}},\quad
    {\bb {m+n}m} = \prod_{d\ge 1} \phi _d^{\fd{m+n}-\fd{m}-\fd{n}}.
  \end{equation*}
  for $m,n\ge 0$.  Here $\floor{r}=\max\{i\in \modZ \ver i\le r\}$ for $r\in \modQ $.
  Set for $a,b,x,y\ge 0$
  \[
  \theta _d(a,b,x,y) = \psi _d(a,x,y)+ \tau _d(x,y)+\tau _d(b,x+y)+\tau _d(a,a+b+x+y),
  \]
  where
  \begin{gather*}
    \psi _d(a,x,y) =\fd{a+x+y}+\fd{a}-\fd{a+x}-\fd{a+y},\\
    \tau _d(x,y) =  \fd{x+y}-\fd{x}-\fd{y}.
  \end{gather*}
  We have
  \begin{equation*}
    \theta (a,b,x,y) = \prod_{d\ge 1}\phi _d^{\theta _d(a,b,x,y)}
  \end{equation*}

  It suffices to prove that if $d\ge 1$ and $a,b,x,y\ge 0$, then
  $\theta _d(a,b,x,y)\ge 0$.  Suppose for contradiction that there are
  integers $d\ge 1$, $a,b,x,y\ge 0$ such that $\theta _d(a,b,x,y)<0$.
  Since $\tau _d(r,s)\ge 0$ for all $r,s\ge 0$, we have
  \begin{equation}
    \label{eq:148}
    \psi _d(a,x,y)<0.
  \end{equation}

  Note that if $a',b',x',y'\ge 0$ are congruent mod $d$ to $a,b,c,d$,
  respectively, then
  \begin{equation*}
    \theta _d(a',b',x',y')=\theta _d(a,b,x,y),\ {}
    \psi _d(a',x',y')=\psi _d(a,x,y),\ {}
    \tau _d(x',y')=\tau _d(x,y).
  \end{equation*}
  For $n\in \modZ $, let $\wt n$ denote the unique integer
  $m\in \{0,1,\ldots,d-1\}$ such that $m\equiv n \pmod d$.  For an
  inequality $P$, let $[P]$ stand for $1$ if $P$ holds and $0$
  otherwise.
  We have
  \begin{equation}
    \label{eq:xi1}
    \begin{split}
      &\psi _d(a,x,y)
      =\psi _d(\ta,\tx,\ty)=\fd{\ta+\tx+\ty}+\fd{\ta}-\fd{\ta+\tx}-\fd{\ta+\ty}\\
      &\quad = [d\le \ta+\tx+\ty]+[2d\le \ta+\tx+\ty]-[d\le \ta+\tx]-[d\le \ta+\ty].
    \end{split}
  \end{equation}
  Since either $d\le \ta+\tx$ or $d\le \ta+\ty$ implies $d\le \ta+\tx+\ty$,
  it follows from \eqref{eq:148} and \eqref{eq:xi1} that
  \begin{gather}
    \label{eq:xi98} \psi _d(a,x,y)=-1,\\
    \label{eq:xi95} d\le \ta+\tx,\  d\le \ta+\ty,\ \ta+\tx+\ty<2d.
  \end{gather}
  By \eqref{eq:xi98} and the assumption $\theta _d(a,b,x,y)<0$,
  we have
  \begin{gather}
    \label{eq:xi5} \tau _d(x,y)=0,\\
    \label{eq:xi6} \tau _d(b,x+y)=0,\\
    \label{eq:xi7} \tau _d(a,a+b+x+y)=0.
  \end{gather}
  By \eqref{eq:xi5}, we have $0=\tau _d(x,y)=\tau _d(\tx,\ty)=[d\le \tx+\ty]$, and
  hence
  \begin{equation}
    \label{eq:xi4} \tx+\ty<d.
  \end{equation}
  Set $z=x+y$.  In view of \eqref{eq:xi4}, we have $\tz=\tx+\ty$.
  By \eqref{eq:xi95}, we have $2d-2\ta\le \tx+\ty<d$,
  and hence
  \begin{gather}
    \label{eq:xi9} 2d-2\ta\le \tz.
  \end{gather}

  By \eqref{eq:xi6}, we have
  $0=\tau _d(b,z)=\tau _d(\tb,\tz)=[d\le \tb+\tz]$, and hence
  \begin{equation}
    \label{eq:xi11} \tb+\tz<d
  \end{equation}
  Set $w=b+z=b+x+y$.  By \eqref{eq:xi11}, we have
  $\tw=\wt{b+z}=\tb+\tz$.  By \eqref{eq:xi9}, we have
  $2d-2\ta\le \tz\le \tz+\tb=\tw$, and hence $2d\le 2\ta+\tw$.
  Therefore
  \[
  \tau _d(a,a+w)=\tau _d(\ta,\ta+\tw)
  =\fd{2\ta+\tw}-\fd{\ta+\tw}-\fd{\ta}=2-\fd{\ta+\tw} - 0 \ge  1,
  \]
  which contradicts to \eqref{eq:xi7}.
  This completes the proof.
\end{proof}

\section{$(\modZ /2)$-gradings in $v$}
\label{sec:q-forms}

The $\modQ (v)$-algebra $U$ admits a ``$\modQ (q)$-form'' $U_q$,
which is the $\modQ (q)$-subalgebra of $U$ generated by the elements
$K$, $K^{-1}$, $vE$, and $F$, and inherits from $U$ a Hopf
$\modQ (q)$-algebra structure.  We have the $(\modZ /2)$-graded
$\modQ (q)$-algebra structure $U=U_q\oplus vU_q$.  We
will call this $(\modZ /2)$-grading the {\em $v$-grading}.

The $\modA $-subalgebras $\UA$, $\modU $, and $\bU$ of $U$ are homogeneous
in the $v$-grading, and hence admits $\modA _q$-forms $\UAq=\UA\cap U_q$,
$\modU _q=\modU \cap U_q$, and $\bU_q=\bU\cap U_q$, respectively.
In particular, $\modU _q$ is the $\modA _q$-subalgebra of $\modU $ generated by
the elements $K$, $K^{-1}$, $e$, and the $v^{-\hf n(n-1)}F^{(n)}$ for
$n\ge 1$, and inherits from $\modU $ a Hopf $\modA _q$-algebra structure.  The
degree $0$ part $\modG _0\modU _q$ of $\modU _q$ is generated by the elements
$K^2$, $K^{-2}$, $e$, and the $v^{-\hf n(n-1)}F^{(n)}K^n$ for $n\ge 1$.

Since the ideals $\modU _n$ are homogeneous in the $v$-grading, the
algebra $\hU$ has a $(\modZ /2)$-graded algebra structure $\hU=\hU_q
\oplus v\hU_q$ with $\hU_q=\plim_n \modU _q/(\modU _n\cap \modU _q)$.
Similarly, $\dU$, $\cU$, and $\tU$ is $(\modZ /2)$-graded with the degree
$0$ part being
\begin{gather*}
  \dU_q = \plim_n \modU _q/(\modU _1\cap \modU _q)^n,\quad
  \cU_q = \plim_n \modU _q/\modU _qe^n\modU _q,\quad
  \tU_q = \operatorname{Im}(\cU_q\rightarrow \hU_q).
\end{gather*}

Since the ideals $\modU _n$, $(\modU _1)^n$, and $\modU ^e_n$ ($n\ge 0$)
are homogeneous in the $v$-grading, the algebra $\hU$ (resp. $\dU$,
$\cU$, $\tU$) has an $\hA_q$- (resp. $\dA_q$-, $\modA _q$-, $\modA _q$-)form,
which we will denote by $\hU_q$ (resp. $\dU_q$, $\cU_q$, $\tU_q$).
These algebras can be regarded as completions of $\modU _q$:
\begin{gather*}
  \hU_q = \plim_n \modU _q/(\modU _n\cap \modU _q),\quad
  \dU_q = \plim_n \modU _q/(\modU _1\cap \modU _q)^n,\\
  \cU_q = \plim_n \modU _q/\modU _qe^n\modU _q,\quad
  \tU_q = \operatorname{Im}(\cU_q\rightarrow \hU_q).
\end{gather*}

\begin{remark}
  Each term of the quasi-$R$-matrix \eqref{eq:90} is contained in
  $\modU _q\otimes _{\modA _q}\modU _q\subset \modU \otimes _\modA \modU $.  (Note that
  $\modU \otimes _\modA \modU =(\modU _q\otimes _{\modA _q}\modU _q)\oplus v(\modU _q\otimes _{\modA _q}\modU _q)$.)
 Therefore we can regard $\Theta $ as an
  element of the completed tensor product
  $\cU_q\check\otimes \cU_q\subset \cU\check\otimes \cU$ of two copies of $\cU_q$, and
  hence as an element of the image
  $\tU_q\tilde\otimes \tU_q=\operatorname{Im}(\cU_q\check\otimes \cU_q\rightarrow \hU_q\hat\otimes \hU_q)$,
  where $\hU_q\hat\otimes \hU_q$ is the completed tensor product of two
  copies of $\hU_q$.
\end{remark}

The centers of $\modU $, $\hU$, $\dU$, and $\tU$ are homogeneous in the
$v$-grading, and we have versions of Theorems \ref{thm:29},
\ref{thm:33}, \ref{thm:34}, \ref{thm:12}, \ref{thm:15} for $\modU _q$ and
its completions.  In particular, we have the following, which we need
in \cite{H:unified}.

\begin{theorem}
  \label{thm:38}
  We have
  \begin{gather*}
    Z(\tU_q) \simeq \plim_n\modA _q[vC]/(\sigma _n),\\
    Z(\modG _0\tU_q)=\modG _0Z(\tU_q) \simeq \plim_n\modA _q[C^2]/(\sigma _n).
  \end{gather*}
\end{theorem}

\end{document}